\documentclass[12pt]{article}
\usepackage{amsmath}
\usepackage{graphicx}
\usepackage{CJK}
\usepackage{multirow}
\usepackage{makecell}

\usepackage{subfigure}
\usepackage[justification=centering]{caption}
\input{amssym.tex}

\def\bc{\begin{center}}       \def\ec{\end{center}}
\def\ba{\begin{array}}        \def\ea{\end{array}}
\def\be{\begin{equation}}     \def\ee{\end{equation}}
\def\bea{\begin{eqnarray}}    \def\eea{\end{eqnarray}}
\def\beaa{\begin{eqnarray*}}  \def\eeaa{\end{eqnarray*}}

\def\mathbb{\Bbb}

\begin{document}
\baselineskip 18pt
\centerline {\bf \large On the number of limit cycles bifurcating from the }
\vskip 0.1 true cm
\centerline {\bf \large  linear center with an algebraic switching curve}

\vskip 0.3 true cm

\centerline{\bf  Jiaxin Wang, Jinping Zhou, Liqin Zhao$^{*}$}
 \centerline{ School of Mathematical Sciences, Beijing Normal University,} \centerline{Laboratory of Mathematics and Complex Systems, Ministry of
Education,} \centerline{Beijing 100875, The People's Republic of China}

\footnotetext[1]{
* Corresponding author.
E-mail:  zhaoliqin@bnu.edu.cn (L. Zhao).}
\vskip 0.2 true cm

\noindent{\bf Abstract}
This paper studies the family of piecewise linear differential systems in the plane with two pieces separated by a switching curve $y=x^{m}$, where $m>1$ is an arbitrary positive. By analysing the first order Melnikov function, we give an upper bound and an lower bound of the maximum number of limit cycles which bifurcate from the period annulus around the origin under polynomial perturbations of degree $n$. The results shows that the degree of switching curves affect the number of limit cycles.


\vskip 0.5 true cm
\centerline{\bf{ $\S$1.} Introduction and the main results}
\vskip 0.5 true cm
Non-smooth differential systems has been widely used in the field of nature, economics [9], nonlinear oscillations [19], and biology [4]. Thus many scholars began to study the dynamical behaviors of this system in recent years. One of the most important problems is to study the existence and number of limit cycles of non-smooth differential systems.

Piecewise smooth differential systems is an important non-smooth differential systems. Many scholars have studied the number of limit cycles for piecewise smooth systems separated by a straight line. In [3,8,10,13,20], these scholars considered some piecewise smooth differential systems which are defined in two zones separated by $x=0$. Some scholars studied a series piecewise smooth differential systems with a switching line $y=0$, see [16,18,21]. In all of these papers, the two most used methods are Melnikov functions and averaging method, which are established in [8,12] and developed in [7,14,15] respectively.

Now, one of the things that comes to mind is what will happen to the existence and number of limit cycles if piecewise smooth differential system with two zones separated by an algebraic curve of degree $n$?  Recently, in [1], a second order averaged functions has been developed and applied to study the number of limit cycles of piecewise linear differential systems with two zones separated by a cubic curve. For the switching curve of degree $n$, the authors [17] considered the crossing limit cycles of a class of discontinuous piecewise linear differential systems formed by two linear differential systems having only centers with a switching curve $y=x^{n}$. The authors in [6] showed that for each $n\in\mathbb{N}$ there exist piecewise linear differential systems separated by an algebraic curve of degree $n$ having $\left[\frac{n}{2}\right]$ hyperbolic limit cycles. In [22], the author studied a piecewise Near-Hamilton systems separated by $y=\pm x^{2}$.

In this paper, motivated by the above analysis, we study the number of limit cycles bifurcating from the linear center with an algebraic switching curve $y=x^{m}$ by the first order Melnikov function, where $m\geq 1$ is an arbitrary positive.

Consider the following perturbed piecewise smooth differential system
$$
\left(
  \begin{array}{c}
    \dot{x} \\
    \dot{y} \\
  \end{array}
\right)
=\left\{
\begin{aligned}
\left(
  \begin{array}{c}
    y+\epsilon p^+(x,y) \\
    -x+\epsilon q^+(x,y) \\
  \end{array}
\right)
,\ y\geq x^m,\\
\left(
  \begin{array}{c}
    y+\epsilon p^-(x,y) \\
    -x+\epsilon q^-(x,y) \\
  \end{array}
\right)
,\ y<x^m, \\
\end{aligned}
\right.\eqno(1.1)_\epsilon$$
where $$p^{\pm}(x,y)=\sum^{n}_{i+j=0}a^{\pm}_{i,j}x^{i}y^{j},~~q^{\pm}(x,y)=\sum^{n}_{i+j=0}b^{\pm}_{i,j}x^{i}y^{j}$$ are any polynomials of degree $n$.
Let $Z(m,n)$ be the upper bound of the number of limit cycles for Hamilton system $(1.1)_\epsilon$(taking into account the multiplicity). Our main results are as follows.
\begin{figure*}[!!ht]
\centering
{\includegraphics[scale=0.3]{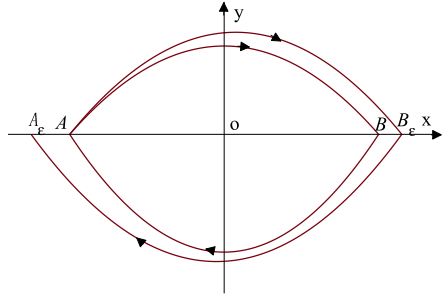}}
\caption{The Poincar\'{e} map related to $y=0$.}
\end{figure*}

\vskip 0.2 true cm
\noindent{\bf Theorem 1.1.}
Suppose that $m=2k+1$($k\in \mathbb{N})$. By using the first order Melnikov function in $\epsilon$, we have

(i)~~For $(m,n)\in D_{1}=:\{(m,n)|0\le n<m-1\}$,
  $$Z(m,n)\leq \frac{1}{2}\left[\frac{n+3}{2}\right]\left[\frac{n+5}{2}\right]+\left[\frac{n+2}{2}\right]\left[\frac{n+4}{2}\right]-2,$$
  $$Z(m,n)\geq \left[\frac{n}{2}\right]\left[\frac{n+6}{2}\right]+\left[\frac{n-1}{2}\right]+2,\qquad \qquad\quad$$
  and the lower bound can be reached by some system $(1.1)_{\epsilon}$.

(ii)~~Let $\delta_n=0$ if $n$ is odd, and $\delta_n=-1$ if $n$ is even. Then we have
$$Z(m,n)\leq \begin{cases}
      \frac{1}{2}\left[\frac{n+3}{2}\right]\left[\frac{n+5}{2}\right]
      +(2k+1)\left[\frac{n}{2}\right]-(k-1)^{2},~\text{if}~(m,n)\in D_{2};\\\\
      (2k+1)\left(2\left[\frac{n}{2}\right]+\delta_n\right)+k(5-3k)+1, ~\text{if}~(m,n)\in D_{3};\end{cases}$$
$$Z(m,n)\geq 2(k+1)\left[\frac{n}{2}\right]+\delta_n-(k-1)^{2}+2, ~~{\text {if $(m,n)\in D_{2}\cup D_{3}$}},  $$
and the lower bound can be reached by some system $(1.1)_{\epsilon}$, where $D_{2}=\{(m,n)|m-1\leq n<2m-1\}$ and $D_{3}=\{(m,n)|n\geq 2m-1\}$.

\vskip 0.2 true cm
\noindent{\bf Theorem 1.2.}
Suppose that $m=2k$($k\in \mathbb{N}$). By using the first order Melnikov function in $\epsilon$, we have

(i)~~ For $(m,n)\in D_{4}=:\{(m,n)|0\leq n<m\}$, we have
  $$Z(m,n)\leq
      4\left[\frac{n}{2}\right]^{2}+\left(6\delta_{n}+11\right)\left[\frac{n}{2}\right]+\frac{1}{2}\delta_{n}\left(5\delta_{n}+17\right)+4,$$
  $$Z(m,n)\geq \frac{1}{2}\left[\frac{n}{2}\right]\left(\left[\frac{n}{2}\right]+3\right)+\frac{1}{2}\left[\frac{n-1}{2}\right]\left(\left[\frac{n-1}{2}\right]+7\right)+3,$$
  and the lower bound can be reached by some system $(1.1)_{\epsilon}$.

(ii)~~ For $(m,n)\in D_{5}\cup D_{6}=:\{(m,n)|m\leq n<2m-2\}\cup\{(m,n)|n\geq 2m-2\}$, we have
  $$Z(m,n)\leq (3k+1)\left[\frac{n-1}{2}\right]+k\left[\frac{n}{2}\right]-k(k-5)-1,$$
  $$Z(m,n)\geq \begin{cases}
  (k+2)\left[\frac{n-1}{2}\right]+k\left[\frac{n}{2}\right]-k(k-3)+1, ~if~(m,n)\in D_{5},\\\\
  3k^{2}+4k-3,~if~n=2m-2,\\\\
  3k^{2}+5k-2,~if~n=2m-1,\\\\
  3k^{2}+6k-2,~if~n=2m.
  \end{cases}$$
  and the lower bound can be reached by some system $(1.1)_{\epsilon}$.

\vskip 0.2 true cm
\noindent{\bf Corollary 1.2.} We have the following results:

(i)~$Z(1,n)=n(n\geq 1)$;

(ii)~$4\left[\frac{n}{2}\right]+2+\delta_n\leq Z(3,n)\leq 6\left[\frac{n}{2}\right]+3(1+\delta_n)$.

 {\noindent}  Particularly, $5\leq Z(3,2)\leq 6$,~$6\leq Z(3,3)\leq 9$ and $9\leq Z(3,4)\leq 12$.

(iii)~$Z(m,1)\geq 2(m\ge 3)$ and  $Z(m,2)\geq 6(m\geq 5)$.

(iv)~$Z(2,1)=3$,~$Z(2,2)=4$ and $3\left(\left[\frac{n}{2}\right]+1\right)+2\delta_{n}\leq Z(2,n)\leq 5\left[\frac{n}{2}\right]+4\delta_{n}+3(n\geq 2)$.

The organizational structure of this paper is as follows. In section 2, we will give some preliminaries. In section 3, we introduce the first order Melnikov function for piecewise smooth near-Hamiltonian system with a switching curve. For Theorem 1.1 and Theorem 1.2, we prove them in section 3 and section 4 respectively.

\vskip 0.2 true cm
\centerline{\bf{ $\S$2}. Preliminaries}
\vskip 0.5 true cm

We first introduce the first order Melnikov function of discontinuous differential systems. Consider the following Hamilton system:
$$
(\dot{x},\ \dot{y})=\begin{cases}
       (H^{+}_{y}(x,y)+\epsilon p^+(x,y),-H^{+}_{x}(x,y)+\epsilon q^+(x,y)),\ \ y>0,\\
       (H^{-}_{y}(x,y)+\epsilon p^-(x,y),-H^{-}_{x}(x,y)+\epsilon q^-(x,y)),\ \ y<0,
\end{cases}
\eqno(2.1)$$
where $0<|\epsilon|\ll 1$, and $p^\pm(x,y)$ and $q^\pm(x,y)$ are polynomials with degree $n$. System $(2.1)$ has two subsystems:
$$
\left\{{\begin{aligned}
\dot{x}&=H^{+}_{y}(x,y)+\epsilon p^{+}(x,y),\\
\dot{y}&=-H^{+}_{x}(x,y)+\epsilon q^{+}(x,y),
\end{aligned}}~~~~~~y>0,
\right.\eqno(2.2)$$
and
$$
\left\{{\begin{aligned}
\dot{x}&=H^{-}_{y}(x,y)+\epsilon p^{-}(x,y),\\
\dot{y}&=-H^{-}_{x}(x,y)+\epsilon q^{-}(x,y),
\end{aligned}}~~~~~~y<0.
\right.\eqno(2.3)$$

We suppose that $(2.1)_{\epsilon=0}$ has a family of periodic orbits around the origin and satisfies the following two assumptions.

\vskip 0.3 true cm

{\bf Assumption (I).} There exist an interval $\Sigma=(\alpha, \beta)$, and two points $A(h)=(a(h),0)$ and $B(h)=(b(h),0)$ such that for $h\in{\Sigma}$
$$
H^{+}(A(h))=H^{+}(B(h))=h,~~H^{-}(A(h))=H^{-}(B(h))=\tilde{h},~~a(h)<b(h).$$

{\bf Assumption (II).} The subsystem $(2.2)_{\epsilon=0}$ has an orbital arc $L_{h}^{+}$ starting from $A(h)$ and ending at $B(h)$ defined by $H^{+}(x,y)=h$ ($y\geq0$). The subsystem $(2.3)_{\epsilon=0}$ has an orbital arc $L_{h}^{-}$ starting from $B(h)$ and ending at $A(h)$ defined by $H^{-}(x,y)=\tilde{h}$($y<0$).

\vskip 0.3 true cm

Under the Assumptions (I) and (II), $(2.1)_{\epsilon=0}$ has a family of non-smooth periodic orbits $L_{h}=L_{h}^{+}\cup L_{h}^{-}(h\in \Sigma)$. For definiteness, we assume that the orbits $L_{h}$ for $h\in \Sigma$ orientate clockwise(see Fig.\,1). The authors [12] established a bifurcation function $F(h,\epsilon)$ for $(2.1)$. Let $F(h,0)=M(h)$. In [8] and [12], the authors obtained the following results.

\begin{figure*}[!!ht]
\centering
{\includegraphics[scale=0.3]{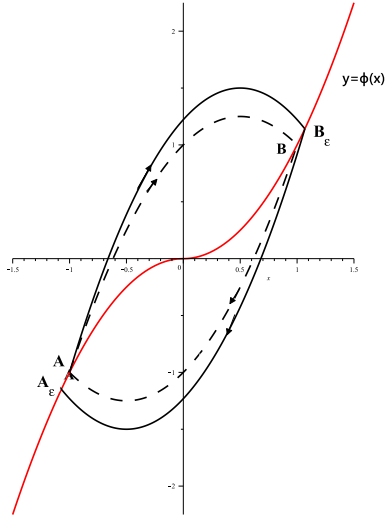}}
\caption{The Poincar\'{e} map related to $y=\phi(x)$.}
\end{figure*}

\vskip 0.3 true cm
\noindent{\bf Lemma 2.1.}([8,12]). Under the assumptions (I) and (II), we have

(i) If $M(h)$ has $k$ zeros in $h$ on the interval $\Sigma$ with each having an odd multiplicity, then $(2.1)$ has at least $k$ limit cycles bifurcating from the period annulus for $0<\left|\epsilon\right|\ll1$.

(ii) If $M(h)$ has at most $k$ zeros in $h$ on the interval $\Sigma$, taking into account the multiplicity, then there exist at most $k$ limit cycles of $(2.1)$ bifurcating from the period annulus.

(iii) The first order Melnikov function $M(h)$ of system $(2.1)$ has the following form
$$
\begin{aligned}
M(h)&=\frac{H_{x}^{+}(A)}{H_{x}^{-}(A)}\left[\frac{H_{x}^{-}(B)}{H_{x}^{+}(B)}\int_{L_{h}^{+}}q^{+}dx-p^{+}dy+\int_{L_{h}^{-}}q^{-}dx-p^{-}dy\right].
\end{aligned}
\eqno(2.4)$$

Further, similar to the [12], if $M(h_{0})=0$ and $M^{'}(h_{0})\neq0$ for some $h_{0}\in{\Sigma}$, then for $|\epsilon|$ small enough system (2.1) has a unique limit cycle near $L_{h_{0}}$. If $M(h)$ has $k$ zeros in $h\in{\Sigma}$ with each having an odd multiplicity, then $(2.1)$ has at least $k$ limit cycles bifurcating from the period annulus for $0<\left|\epsilon\right|\ll1$. If $M(h)$ has at most $k$ zeros in $h\in{\Sigma}$, taking into account the multiplicity, then there exist at most $k$ limit cycles of $(2.1)$ bifurcating from the period annulus.

\vskip 0.2 true cm
\noindent{\bf Lemma 2.2.}[2] Consider $p+1$ linearly independent analytical functions $f_i:U\rightarrow \mathbb{R}, i=0,1,...,p$, where $U \in \mathbb{R}$ is an interval. Suppose that there exists $j \in \{0,,1,...,p\}$ such that $f_{j}$ has constant sign. Then there exists $p+1$ constants $C_{i} ,i=0,1,...,p$ such that $f(x)=\sum _{i=0}^{p}C_{i}f_{i}(x)$ has at least $p$ simple zeros in $U$.

\vskip 0.2 true cm

\noindent{\bf Definition 2.3.}[5] Let $p_0(x),p_1(x),...,p_{n-1}(x)$ be analytic functions on an open interval $J\subset\mathbb{R}$. The ordered set $\left(p_0(x),p_1(x),...,p_{n-1}(x)\right)$ is said to be an ECT-system on $J$ if, for all $k=1,2,...,n$, any nontrivial linear combination
$$\alpha_0p_0(x)+\alpha_1p_1(x)+...+\alpha_{k-1}p_{k-1}(x)$$
has at most $k-1$ isolated zeros on $J$ counted with multiplicities.

\vskip 0.2 true cm
\noindent{\bf Lemma 2.4.}[5] The ordered set $(p_0(x),p_1(x),...,p_{n-1}(x))$ is an ECT-system on $J$ if and only if, for each $k=1,2,...,n,$
$$W(p_0,p_1,...,p_{k-1})\neq0,$$ for all $x\in J,$
where $W(p_0,p_1,...,p_{k-1})$ is the Wronskian of functions $p_0(x),p_1(x),...,p_{k-1}(x).$

\vskip 0.5 true cm

\centerline{\bf{ $\S$3}. Proof of Theorem 1.1.}
\vskip 0.5 true cm

We first introduce the first order Melnikov function for piecewise smooth near-Hamiltonian system with a switching curve defined by $y=\phi(x)$. Consider the following system
$$
\left(
  \begin{array}{c}
    \dot{x} \\
    \dot{y} \\
  \end{array}
\right)
=\left\{
\begin{aligned}
\left(
  \begin{array}{c}
    H_{y}^{+}(x,y)+\epsilon f^+(x,y) \\
    -H_{x}^{+}(x,y)+\epsilon g^+(x,y) \\
  \end{array}
\right)
,\ y\geq \phi(x),\\
\left(
  \begin{array}{c}
    H_{y}^{-}(x,y)+\epsilon f^-(x,y) \\
    -H_{x}^{-}(x,y)+\epsilon g^-(x,y) \\
  \end{array}
\right)
,\ y<\phi(x)\\
\end{aligned}
\right.\eqno(3.1)
$$
where $H^{\pm}$, $f^{\pm}$, $g^{\pm}$ and $\phi(x)$ are all $C^{\infty}$ functions satisfying $\phi(0)=0$, $\epsilon\geq0$ is a small parameter. System $(3.1)$ has two subsystems:
$$
\left\{{\begin{aligned}
\dot{x}&=H^{+}_{y}(x,y)+\epsilon f^{+}(x,y),\\
\dot{y}&=-H^{+}_{x}(x,y)+\epsilon g^{+}(x,y),
\end{aligned}}~~~~~~y\geq \phi(x),
\right.\eqno(3.2)$$
and
$$
\left\{{\begin{aligned}
\dot{x}&=H^{-}_{y}(x,y)+\epsilon f^{-}(x,y),\\
\dot{y}&=-H^{-}_{x}(x,y)+\epsilon g^{-}(x,y),
\end{aligned}}~~~~~~y<\phi(x).
\right.\eqno(3.3)$$

We suppose that $(3.1)_{\epsilon=0}$ has a family of periodic orbits around the origin and satisfies the following assumptions.

\noindent{\bf (A1)}: There exists an open interval $\Sigma$ such that for each $h\in\Sigma$, there are two points $A(h)$ and $B(h)$ on the curve $y=\phi(x)$ with $A(h)=(a(h),\phi(a(h)))$, $B(h)=(b(h),\phi(b(h)))$ and satisfying
$$
H^{+}(A(h))=H^{+}(B(h))=h,~~H^{-}(A(h))=H^{-}(B(h)),~~a(h)<0<b(h).
$$

\noindent{\bf (A2)}: The subsystem $(3.2)_{\epsilon=0}$ has an orbital arc $L_{h}^{+}$ starting from $A(h)$ and ending at $B(h)$ defined by $H^{+}(x,y)=h$ ($y\geq \phi(x)$). The subsystem $(3.3)_{\epsilon=0}$ has an orbital arc $L_{h}^{-}$ starting from $B(h)$ and ending at $A(h)$ defined by $H^{-}(x,y)=\tilde{h}$($y<\phi(x)$).

\noindent{\bf (A3)}: Curve $L_{h}^{\pm}$, $h\in\Sigma$ are not tangent to curve $y=\phi(x)$ at points $A(h)$ and $B(h)$. In other words, for each $h\in\Sigma$,
$$H_{x}^{\pm}(x,y)+H_{y}^{\pm}(x,y)\phi^{'}(x)\neq0$$
at points $A(h)$ and $B(h)$.

Under the Assumptions ${\bf (A1)}$, ${\bf (A2)}$ and ${\bf (A3)}$, $(2.1)_{\epsilon=0}$ has a family of non-smooth periodic orbits $L_{h}=L_{h}^{+}\cup L_{h}^{-}(h\in\Sigma)$. For definiteness, we assume that the orbits $L_{h}$ for $h\in\Sigma$ orientate clockwise(see Fig.\,2). We have the following Lemma.

\vskip 0.2 true cm
\noindent{\bf Lemma 3.1.} The first order Melnikov function $M(h)$ of system (3.1) can be expressed as
$$M(h)=\int_{L_{h}^{+}}g^{+}dx-f^{+}dy+\frac{H_{x}^{+}(A)+H_{y}^{+}(A)\phi^{'}(a(h))}{H_{x}^{-}(A)+H_{y}^{-}(A)\phi^{'}(a(h))}\int_{L_{h}^{-}}g^{-}dx-f^{-}dy.
\eqno(3.4)$$

\vskip 0.2 true cm
\noindent{\bf Proof}. Let us make the following transformation:
$$x=x,~~z=y-\phi(x).$$
Then system (3.1) is reduced to
$$
\left(
  \begin{array}{c}
    \dot{x} \\
    \dot{z} \\
  \end{array}
\right)
=\left\{
\begin{aligned}
\left(
  \begin{array}{c}
    {\widetilde H}_{z}^{+}(x,z)+\epsilon p^+(x,z) \\
    -{\widetilde H}_{x}^{+}(x,z)+\epsilon q^+(x,z) \\
  \end{array}
\right)
,\ z\geq 0,\\
\left(
  \begin{array}{c}
    {\widetilde H}_{z}^{-}(x,z)+\epsilon p^-(x,z) \\
    -{\widetilde H}_{x}^{-}(x,z)+\epsilon q^-(x,z) \\
  \end{array}
\right)
,\ z<0, \\
\end{aligned}
\right.\eqno(3.5)
$$
where $${\widetilde H}^{\pm}(x,z)=H^{\pm}(x,z+\phi(x)), ~{\widetilde H}^{\pm}_{z}(x,z)=H^{\pm}_{y}(x,z+\phi(x)),$$  $${\widetilde H}^{\pm}_{x}(x,z)=H_{x}^{\pm}(x,z+\phi(x))+H_{y}^{\pm}(x,z+\phi(x))\phi^{'}(x),
\qquad \quad $$ and
$$~ p^{\pm}(x,z)=f^{\pm}(x,z+\phi(x)),\qquad \qquad \qquad$$
$$\qquad \qquad q^{\pm}(x,z)=g^{\pm}(x,z+\phi(x))-\phi^{'}(x)f^{\pm}(x,z+\phi(x)).$$ Denote ${\widetilde A}(h)=(a(h),0)$. By  Lemma 2.1, we have
$$\begin{aligned}
{ M}(h)&=\int_{{\widetilde L}_{h}^{+}}q^{+}(x,z)dx-p^{+}(x,z)dz+\frac{{\widetilde H}^{+}_{x}({\widetilde A})}{{\widetilde H}^{-}_{x}({\widetilde A})}\int_{{\widetilde L}_{h}^{-}}q^{-}(x,z)dx-p^{-}(x,z)dz\\
&=\int_{{\widetilde L}_{h}^{+}}\left[g^{+}(x,z+\phi(x))-\phi^{'}(x)f^{+}(x,z+\phi(x))\right]dx-f^{+}(x,z+\phi(x))d(y-\phi(x))\\
&+\frac{{\widetilde H}^{+}_{x}({\widetilde A})}{{\widetilde H}^{-}_{x}({\widetilde A})}\int_{{\widetilde L}_{h}^{-}}\left[g^{-}(x,z+\phi(x))-\phi^{'}(x)f^{-}(x,z+\phi(x))\right]dx\\
&-f^{-}(x,z+\phi(x))d(y-\phi(x))\\
&=\int_{L_{h}^{+}}g^{+}dx
-f^{+}dy+\frac{H_{x}^{+}(A)+H_{y}^{+}(A)\phi^{'}(a(h))}
{H_{x}^{-}(A)+H_{y}^{-}(A)\phi^{'}(a(h))}\int_{L_{h}^{-}}g^{-}dx-f^{-}dy.
\end{aligned}$$
 This ends the proof. $\diamondsuit$
\vskip 0.3 true cm

Next, we will obtain the algebraic structure of $M(h)$ for  system $(1.1)_{\epsilon}$. For $h\in(0,+\infty)$ and $i,j\in\mathbb{N}$, we denote
$$J_{i,j}(h)=\int_{L_{h}^{+}}x^{i}y^{j}dx,~~~I_{i,j}(h)=\int_{L_{h}^{-}}x^{i}y^{j}dx,$$
where
$$L_{h}^{\pm}=\{(x,y)|H(x,y)={h}/{2},y\geq x^{m}(y\leq x^{m})\}.$$
By Lemma 3.1, we have
$$
\begin{aligned}
M(h)&=\int_{L_{h}^{+}}q^{+}dx-p^{+}dy+\int_{L_{h}^{-}}q^{-}dx-p^{-}dy\\
    &=\sum^{n}_{i+j=0}\left(\int_{L_{h}^{+}}b^{+}_{i,j}x^{i}y^{j}dx-a^{+}_{i,j}x^{i}y^{j}dy+\int_{L_{h}^{-}}b^{-}_{i,j}
    x^{i}y^{j}dx-a^{-}_{i,j}x^{i}y^{j}dy\right).
\end{aligned}
\eqno(3.6)
$$
Suppose that the orbit $L_{h}^{+}(L_{h}^{-})$ intersects the curve $y=x^m$ at points $A(-u(h),(-u(h))^{m})$ and $B(u(h),u(h)^{m})$.

\vskip 0.3 true cm

\noindent{\bf Lemma 3.2.} The first order Melnikov function $M(h)$ can be written as
$$
M(h)=\sum^{n}_{i+j=0}\rho_{i,j}^{+}J_{i,j}(h)+\sum^{n}_{i+j=0}\rho_{i,j}^{-}I_{i,j}(h)+\Phi(u(h)),\eqno(3.7)$$
where $\rho_{i,j}^{\pm}$ are arbitrary constants, and $\Phi(u)$ is a polynomial of $u$ with degree no more than $m(n+1)$.

\vskip 0.2 true cm
\noindent{\bf Proof.}
Using the Green's Formula, we have
$$
\begin{aligned}
\int_{L_{h}^{+}}x^{i}y^{j}dy&=\int_{L_{h}^{+}\cup\widehat{BOA}}x^{i}y^{j}dy-\int_{\widehat{BOA}}x^{i}y^{j}dy\\
&=-i\iint_{int(L_{h}^{+}\cup\widehat{BOA})}x^{i-1}y^{j}dxdy-\frac{m((-1)^{i+mj+m}-1)}{i+mj+m}u(h)^{i+mj+m},\\
\int_{L_{h}^{+}}x^{i}y^{j}dx&=\int_{L_{h}^{+}\cup\widehat{BOA}}x^{i}y^{j}dx-\int_{\widehat{BOA}}x^{i}y^{j}dx\\
&=j\iint_{int(L_{h}^{+}\cup{\widehat{BOA}})}x^{i}y^{j-1}dxdy-\frac{(-1)^{i+mj+1}-1}{i+mj+1}u(h)^{i+mj+1},
\end{aligned}
$$
which imply that
$$
\int_{L_{h}^{+}}x^{i}y^{j}dy=-\frac{i}{j+1}\int_{L_{h}^{+}}x^{i-1}y^{j+1}dx-\frac{(-1)^{i+mj+m}-1}{j+1}u(h)^{i+mj+m}.\eqno(3.8)
$$
In the similar way, we have
$$
\int_{L_{h}^{-}}x^{i}y^{j}dy=-\frac{i}{j+1}\int_{L_{h}^{-}}x^{i-1}y^{j+1}dx-\frac{1-(-1)^{i+mj+m}}{j+1}u(h)^{i+mj+m}.\eqno(3.9)
$$

From (3.6), (3.8) and (3.9), we can obtain
$$
\begin{aligned}
M(h)&=\sum\limits^{n}_{i+j=0}\left(\int_{L_{h}^{+}}\left(b^{+}_{i,j}x^{i}y^{j}+\frac{i}{j+1}a^{+}_{i,j}x^{i-1}y^{j+1}\right)dx\right.\\
&\left.+\int_{L_{h}^{-}}\left(
b^{-}_{i,j}x^{i}y^{j}+\frac{i}{j+1}a^{-}_{i,j}x^{i-1}y^{j+1}\right)dx\right)\\
&+\sum^{n}_{i+j=0}(a^{+}_{i,j}-a^{-}_{i,j})\frac{(-1)^{i+mj+m}-1}{j+1}u(h)^{i+mj+m}\\
:&=\sum^{n}_{i+j=0}\rho^{+}_{i,j}J_{i,j}(h)+\sum^{n}_{i+j=0}\rho^{-}_{i,j}I_{i,j}(h)+\Phi(u(h)),
\end{aligned}
$$
where $\Phi(u(h))=\sum\limits^{n}_{i+j=0}(a^{+}_{i,j}-a^{-}_{i,j})\frac{(-1)^{i+mj+m}-1}{j+1}u(h)^{i+mj+m}$, $\rho^{\pm}_{i,j}=b^{\pm}_{i,j}+\frac{i+1}{j}a^{\pm}_{i+1,j-1}(j\geq1)$ and $\rho^{\pm}_{i,0}=b^{\pm}_{i,0}$. This ends the proof. $\diamondsuit$

From now on, we will consider the case for $m=2k+1$. Denote $J_{i,j}(h),~I_{i,j}(h)$ and $u(h)$ as $J_{i,j},~I_{i,j}$ and $u$.
\vskip 0.2 true cm
\noindent{\bf Lemma 3.3.}
For $h\in(0,+\infty)$, investigating $J_{i,j}$ and $I_{i,j}(i+j=n)$, we have
(i)~~For $i,l,s\geq0$, we have
$$\begin{aligned}
J_{i,2s}=-I_{i,2s},~J_{2l,2s+1}=I_{2l,2s+1},~J_{2l+1,2s+1}=-I_{2l+1,2s+1}.
\end{aligned}$$
(ii)~~ If $n=2l$, then for $d\geq 0$, we have
  $$\begin{aligned}
  J_{2l-2d,2d}&=\tau_{l,d}^{0}h^{l}J_{0,0}+\sum_{i=1}^{d}\tau_{l,d}^{i}h^{i-1}u^{2l+2d(m-1)+1-2m(i-1)}\\
  &+\sum_{i=1}^{l-d}\tau_{l,d}^{i+d}h^{l-i}u^{2m+1+2(i-1)},
  \end{aligned}\eqno(3.10)$$
  $$\begin{aligned}
  J_{2l-2d-1,2d+1}&=\mu_{l,d}^{0}h^{l-1}J_{1,1}+\sum_{i=1}^{d}\mu_{l,d}^{i}h^{i-1}u^{2l+2d(m-1)+m-2m(i-1)}\\
  &+\sum_{i=2}^{l-d}\mu_{l,d}^{i+d}h^{l-i}u^{3m+2+2(i-2)}.
  \end{aligned}\eqno(3.11)$$

(iii)~~If $n=2l+1$, then for $d\geq 0$, we have $$J_{2l+1-2d,2d}=0,~
  J_{2l-2d,2d+1}=\chi_{l,d}h^{l}J_{0,1},\eqno(3.12)$$
where $\tau_{l,d}^{i},~\mu_{l,d}^{i}$ and $\chi_{l,d}$ are arbitrary constants.
\vskip 0.2 true cm
\noindent{\bf Proof.}
According to the definition of $L_{h}^{\pm}$, it follows that
$$
J_{i,2s}=\int_{L_{h}^{+}}x^{i}y^{2s}dx=\int^{u}_{-u}x^{i}(h-x^{2})^{s}dx=-\int^{-u}_{u}x^{i}(h-x^{2})^{s}dx=-I_{i,2s},
$$
$$\begin{aligned}
J_{2r,2s+1}&=\int_{L_{h}^{+}}x^{2r}y^{2s+1}dx=-\int_{-u}^{-\sqrt{h}}x^{2r}(\sqrt{h-x^{2}})^{2s+1}dx+\int_{-\sqrt{h}}^{u}x^{2r}(\sqrt{h-x^{2}})^{2s+1}dx\\
&=\int_{u}^{\sqrt{h}}x^{2r}(\sqrt{h-x^{2}})^{2s+1}dx+\int_{-u}^{\sqrt{h}}x^{2r}(\sqrt{h-x^{2}})^{2s+1}dx=I_{2r,2s+1}.
\end{aligned}$$
It is similar with $J_{2r+1,2s+1}=-I_{2r+1,2s+1}$.
Differentiating $H(x,y)=x^{2}+y^{2}=\frac{h}{2}$ with respect to $x$, we obtain
$$
x+y\frac{\partial y}{\partial x}=0.
\eqno(3.13)$$
Multiplying $H(x,y)=\frac{h}{2}$ and (3.13) by $x^{i}y^{j}dx$ and $x^{i+1}y^{j}dx$ respectively and integrating over $L_{h}^{+}$, noting (3.8) we have
$$
J_{i+2,j}+J_{i,j+2}=hJ_{i,j}
\eqno(3.14)$$
$$
J_{i+2,j}-\frac{i+1}{j+2}J_{i,j+2}-\frac{(-1)^{i+mj+2m+1}-1}{j+2}u^{i+mj+2m+1}=0.
\eqno(3.15)$$
Elementary manipulations reduce Eps. (3.14) and (3.15) to
$$
J_{i,j}=\frac{j}{i+j+1}\left(hJ_{i,j-2}-\frac{(-1)^{i+mj+1}-1}{j}u^{i+mj+1}\right),\eqno(3.16)$$
$$
J_{i,j}=\frac{j+2}{i+j+1}\left(\frac{i-1}{j+2}hJ_{i-2,j}+\frac{(-1)^{i+mj+2m-1}-1}{j+2}u^{i+mj+2m-1}\right).
\eqno(3.17)$$
We will prove the conclusion by induction on $n$. Without loss of generality, we only prove (3.10), the rest can be shown in a similar way. It is obvious that $J_{1,0}=0$. Then when $l=1,2,3$, (3.16) and (3.17) give
$$J_{2,0}=\frac{2}{3}hJ_{0,0}-\frac{2}{3}u^{2m+1},$$
$$J_{0,2}=\frac{2}{3}hJ_{0,0}+\frac{2}{3}u^{2m+1},$$
$$J_{4,0}=h^{2}J_{0,0}-\frac{2}{3}hu^{2m+1}-\frac{2}{5}u^{2m+3},$$
$$J_{3,1}=\frac{4}{5}hJ_{1,1}-\frac{2}{5}u^{3m+2},$$
$$J_{1,3}=\frac{3}{5}hJ_{1,1}+\frac{2}{5}u^{3m+2},$$
$$J_{2,2}=\frac{2}{5}h^{2}J_{0,0}-\frac{4}{15}hu^{2m+1}+\frac{2}{5}u^{2m+3},$$
$$J_{0,4}=\frac{8}{15}h^{2}J_{0,0}+\frac{8}{15}hu^{2m+1}+\frac{2}{5}u^{4m+1},$$
$$J_{6,0}=h^{3}J_{0,0}-\frac{2}{3}h^{2}u^{2m+1}-\frac{2}{5}hu^{2m+3}
-\frac{2}{7}u^{2m+5},$$
$$J_{5,1}=\frac{24}{35}h^{2}J_{1,1}-\frac{12}{35}hu{3m+2}-\frac{1}{7}u^{3m+4},$$
$$J_{4,2}=\frac{2}{7}h^{3}J_{0,0}-\frac{4}{21}h^{2}u^{2m+1}-\frac{4}{35}h
u^{2m+3}+\frac{2}{7}u^{2m+5},$$
$$J_{3,3}=\frac{12}{35}h^{2}J_{1,1}-\frac{6}{35}hu^{3m+2}+\frac{2}{7}u^{3m+4},$$
$$J_{2,4}=\frac{8}{35}h^{3}J_{0,0}-\frac{16}{105}h^{2}u^{2m+1}
+\frac{8}{35}hu^{2m+3}+\frac{2}{7}u^{4m+3},$$
$$J_{1,5}=\frac{3}{7}h^{2}J_{1,1}+\frac{2}{7}hu^{3m+2}+\frac{2}{7}u^{5m+2},$$
$$J_{0,6}=\frac{4}{35}h^{3}J_{0,0}+\frac{4}{35}h^{2}u^{2m+1}
+\frac{12}{35}hu^{4m+1}+\frac{2}{7}u^{6m+1},$$
which yield the conclusion for $=1,2,3$. Suppose that the result holds for $l\leq k-1$. Then for $l=k$, taking $d=0$ in (3.17), $d=1,...,k-1,k$ in (3.16) respectively, we can obtain that
$$
\left(\begin{matrix}
          J_{2k,0}\\
          J_{2k-2,2}\\
          \vdots\\
          J_{2,2k-2}\\
          J_{0,2k}
          \end{matrix}\right)\ \
=\left(\begin{matrix}
        hJ_{2k-2,0}-u^{2k+2m-1}\\
        \frac{2}{2k+1}hJ_{2k-2,0}+u^{2k+2m-1}\\
        \vdots\\
        \frac{2k-2}{2k+1}hJ_{2,2k-4}+\frac{1}{k-1}u^{(2k-2)m+3}\\
        \frac{2k}{2k+1}hJ_{0,2k-2}+\frac{1}{k}u^{2km+1}
        \end{matrix}\right).\eqno(3.18)$$
By inductive hypothesis, we have for $i+j=2k$ and $d\geq 0$,
$$
J_{2k-2d,2d}=\tau_{k,d}^{0}h^{k}J_{0,0}+\sum_{i=1}^{d}\tau_{k,d}^{i}h^{i-1}u^{2k+2d(m-1)+1-2m(i-1)}+\sum_{i=1}^{k-d}\tau_{k,d}^{i+d}h^{k-i}u^{2m+1+2(i-1)}.
$$
This completes the proof of Lemma 3.2. $\diamondsuit$

Let $\rho_{i,j}=\rho_{i,j}^{+}-\rho_{i,j}^{-}$, $\gamma_{i,j}=a_{i,j}^{+}-a_{i,j}^{-}$ and $\zeta_{i,j}=\rho_{i,j}^{+}+\rho_{i,j}^{-}$. By Lemma 3.2 and Lemma 3.3, we can get
$$\begin{aligned}
M(h)&=\sum_{l=0}^{\left[\frac{n}{2}\right]}\left(\sum_{i+j=2l}(\rho_{i,j}^{+}-\rho_{i,j}^{-})J_{i,j}+\sum_{i+j=2l}(a_{i,j}^{+}-a_{i,j}^{-})\frac{(-1)^{i+mj+m}-1}{j+1}u^{i+mj+m}\right)\\
&+\sum_{l=0}^{\left[\frac{n-1}{2}\right]}\left(\sum_{k+s=2l+1}(\rho_{k,s}^{+}+\rho_{k,s}^{-})J_{k,s}+\sum_{k+s=2l+1}(a_{k,s}^{+}-a_{k,s}^{-})\frac{(-1)^{k+ms+m}-1}{s+1}u^{k+ms+m}\right)\\
&=:\sum_{l=0}^{\left[\frac{n}{2}\right]}\sum_{i+j=2l}\left(\rho_{i,j}J_{i,j}+\gamma_{i,j}\frac{(-1)^{i+mj+m}-1}{j+1}u^{i+mj+m}\right)+\sum_{l=0}^{\left[\frac{n-1}{2}\right]}
\sum_{k+s=2l+1}\zeta_{k,s}J_{k,s}.
\end{aligned}$$
It is easy to show that the coefficients $\rho_{i,j}$, $\gamma_{i,j}(i+j=2l)$ and $\zeta_{k,s}(k+s=2l+1)$ are independent. In fact, let
$${\bf M}:=\frac{\partial{(\rho_{0,0},\rho_{2,0},\rho_{1,1},...,\rho_{0,2\left[\frac{n}{2}\right]},\gamma_{0,0},...,\gamma_{0,2\left[\frac{n}{2}\right]}
,\zeta_{1,0},...,\zeta_{0,2\left[\frac{n-1}{2}\right]+1})}}{\partial{(b_{0,0}^{+},b_{2,0}^{+},b_{1,1}^{+},...,b_{0,2\left[\frac{n}{2}\right]}^{+},a_{0,0}^{+},...,
a_{0,2\left[\frac{n}{2}\right]}^{+},b_{1,0}^{+},...,b_{0,2\left[\frac{n-1}{2}\right]+1}^{+})}}.$$
 We have ${\rm det}~{\bf M}=1$ by the expression of $\rho_{i,j}$ and  $\gamma_{i,j}$ and $\zeta_{k,s}$. This ends the proof.

\vskip 0.2 true cm

\noindent{\bf Lemma 3.4.} Let $h=u^2+u^{2m}$. Then the number of zeros of  $M(h)$ is equal to the number of zeros of $M(u)$ in $u\in(0,+\infty)$. For $(m,n)\in D_{1}$, $M(u)$ can be written as
  $$M(u)=\sum_{p=0}^{\left[\frac{n}{2}\right]}A_{p}u^{2p+1}+\sum_{p=0}^{\left[\frac{n}{2}\right]}\sum_{k=0}^{2p}B_{p,k}u^{2p_{0}+1}+\sum_{l=0}^{\left[\frac{n-1}{2}\right]}K_{l}(u^2+u^{2m})^{l+1},\eqno(3.19)$$
  where
  $$\begin{aligned}
  p_{0}&=p+\frac{1}{2}(m-1)(k+1),\\
  A_{p}&=\sum\limits_{\substack{l+(m-1)j=p\\j\leq l}}\sum\limits_{d=j}^{l}\omega_{l,d}^{j}\rho_{2l-2d,2d}+\sum\limits_{\substack{l+(s+\frac{3}{2})(m-1)=p\\s\leq l-1}}\sum\limits_{d=0}^{l-1}\vartheta_{l,d}^{s}\rho_{2l-2d-1,2d+1},\\
B_{p,k}&=\sum\limits_{\substack{l+j(m-1)=p_{0}\\j\leq l}}\sum\limits_{d=j}^{l}\omega_{l,d}^{j}\rho_{2l-2d,2d}+\sum\limits_{\substack{l+(s+\frac{3}{2})(m-1)=p_{0}\\s\leq {l-1}}}\sum\limits_{d=0}^{l-1}\vartheta_{l,d}^{s}\rho_{2l-2d-1,2d+1}\\
&-\sum\limits_{i+mj+m=2p_{0}+1}\frac{2}{j+1}\gamma_{i,j},\\
K_{l}&=\sum\limits_{k=0}^{l}\chi_{l,k}\zeta_{2l-2k,2k+1},
  \end{aligned}$$
and $\omega_{l,d}^{j}$, $\vartheta_{l,d}^{s}$ and $\chi_{l,k}$ are constants with $\omega_{l,0}^{0}=\frac{2}{2l+1}$ and $\chi_{l,0}\neq 0$.

\vskip 0.2 true cm
\noindent{\bf Proof.}
 By direct computations, we can obtain
$$J_{0,0}=2u,~~J_{0,1}=\frac{\pi}{2}(u^2+u^{2m}),~~J_{1,1}=-\frac{2}{3}u^{3m}.$$
Noticing $h=u^{2}+u^{2m}$ and  (3.10)--(3.12), for  $i+j=2l$ and  $d\geq 0$,  we have
  $$\begin{aligned}
  J_{2l-2d,2d}&=\tau_{l,d}^{0}(u^{2}+u^{2m})^{l}J_{0,0}+\sum_{i=1}^{d}\tau_{l,d}^{i}(u^{2}+u^{2m})^{i-1}u^{2l+2d(m-1)+1-2m(i-1)}\\
  &+\sum_{i=1}^{l-d}\tau_{l,d}^{i+d}(u^{2}+u^{2m})^{l-i}u^{2m+1+2(i-1)}\\
  &=\sum_{j=0}^{l}2\tau_{l,d}^{0}\binom l j u^{2l+2j(m-1)+1}+\sum_{i=1}^{d}\sum_{j=0}^{i-1}\tau_{l,d}^{i}\binom {i-1} j u^{2(i-1-j-d)(1-m)+2l+1}\\
  &+\sum_{i=1}^{l-d}\sum_{j=0}^{l-i}\tau_{l,d}^{i+d}\binom {l-i} j u^{2(l-1)+2j(m-1)+2m+1}\\
  &=\sum_{j=0}^{l}\omega_{l,d}^{j}u^{2l+2j(m-1)+1},
  \end{aligned}
  \eqno(3.20)$$
  $$\begin{aligned}
  J_{2l-2d-1,2d+1}&=\mu_{l,d}^{0}(u^{2}+u^{2m})^{l-1}J_{1,1}+\sum_{i=1}^{d}\mu_{l,d}^{i}(u^{2}+u^{2m})^{i-1}u^{2l+2d(m-1)+m-2m(i-1)}\\
  &+\sum_{i=2}^{l-d}\mu_{l,d}^{i+d}(u^{2}+u^{2m})^{l-i}u^{3m+2+2(i-1)}\\
  &=\sum_{j=0}^{l-1}-\frac{2}{3}\mu_{l,d}^{0}\binom {l-1} j u^{2(l-1)+2j(m-1)+3m}\\
  &+\sum_{i=1}^{d}\sum_{j=0}^{i-1}\mu_{l,d}^{i}\binom {i-1} j u^{2(i-1-j-d)(1-m)+2l+m}\\
  &+\sum_{i=2}^{l-d}\sum_{j=0}^{l-i}\mu_{l,d}^{i+d}\binom {l-i} j u^{2(l-2)+2j(m-1)+3m+2}\quad\quad\quad\quad\quad\quad\quad\quad\quad\quad\quad\\
  &=\sum_{s=0}^{l-1}\vartheta_{l,d}^{s}u^{2l+2s(m-1)+3m-2}.
  \end{aligned}
  \eqno(3.21)$$
For $i+j=2l+1$, $d\geq 0$ we have
  $$J_{2l-2d,2d+1}=\chi_{l,d}(u^{2}+u^{2m})^{l+1},\eqno(3.22)$$
  $$J_{2l-2d+1,2d}=0,\eqno(3.23)$$
where $\omega_{l,d}^{j}$, $\vartheta_{l,d}^{s}$ and $\chi_{l,d}$ are constants.
Therefore, we can obtain
$$\begin{aligned}
M(h)&=\sum_{l=0}^{\left[\frac{n}{2}\right]}\sum_{i+j=2l}\left(\rho_{i,j}J_{i,j}+\gamma_{i,j}\frac{(-1)^{i+mj+m}-1}{j+1}u^{i+mj+m}\right)+\sum_{l=0}^{\left[\frac{n-1}{2}\right]}\sum_{i+j=2l+1}\zeta_{i,j}J_{i,j}\\
&=\sum_{l=0}^{\left[\frac{n}{2}\right]}\sum_{d=0}^{l}\rho_{2l-2d,2d}J_{2l-2d,2d}+\sum_{l=0}^{\left[\frac{n}{2}\right]}\sum_{d=0}^{l-1}\rho_{2l-2d-1,2d+1}J_{2l-2d-1,2d+1}\\
&+\sum_{l=0}^{\left[\frac{n-1}{2}\right]}\sum_{d=0}^{2l+1}\zeta_{2l-2d,2d+1}J_{2l-2d,2d+1}+\sum_{l=0}^{\left[\frac{n}{2}\right]}\sum_{i+j=2l}\gamma_{i,j}\frac{(-1)^{i+mj+m}-1}{j+1}u^{i+mj+m}\\
&=\sum_{l=0}^{\left[\frac{n}{2}\right]}\sum_{d=j}^{l}\sum_{j=0}^{l}\omega_{l,d}^{j}\rho_{2l-2d,2d}u^{2l+2j(m-1)+1}+\sum_{l=0}^{\left[\frac{n}{2}\right]}\sum_{d=0}^{l-1}\sum_{j=1}^{l}\vartheta_{l,d}^{s}\rho_{2l-2d-1,2d+1}u^{2l+2j(m-1)+m}\\
&+\sum_{l=0}^{\left[\frac{n-1}{2}\right]}\sum_{d=0}^{2l+1}\chi_{l,d}\zeta_{2l-2d,2d+1}(u^{2}+u^{2m})^{l+1}+\sum_{l=0}^{\left[\frac{n}{2}\right]}\sum_{i+j=2l}\gamma_{i,j}\frac{(-1)^{i+mj+m}-1}{j+1}u^{i+mj+m}.
\end{aligned}
$$
By direct analysises and rearrangement, we can obtain for $(m,n)\in D_{1}$,
$$\begin{aligned}
M(u)&=\sum_{p=0}^{\left[\frac{n}{2}\right]}\left(\sum\limits_{\substack{l+(m-1)j=p\\j\leq l}}\sum\limits_{d=j}^{l}\omega_{l,d}^{j}\rho_{2l-2d,2d}+\sum\limits_{\substack{l+(s+\frac{3}{2})(m-1)=p\\s\leq l-1}}\sum\limits_{d=0}^{l-1}\vartheta_{l,d}^{s}\rho_{2l-2d-1,2d+1}\right)u^{2p+1}\\
&+\sum_{p=0}^{\left[\frac{n}{2}\right]}\sum_{k=0}^{2p}\left(\sum\limits_{\substack{l+j(m-1)=p_{0}\\j\leq l}}\sum\limits_{d=j}^{l}\omega_{l,d}^{j}\rho_{2l-2d,2d}+\sum\limits_{\substack{l+(s+\frac{3}{2})(m-1)=p_{0}\\s\leq {l-1}}}\sum\limits_{d=0}^{l-1}\vartheta_{l,d}^{s}\rho_{2l-2d-1,2d+1}\right.\\
&\left.+\sum\limits_{i+mj+m=2p_{0}+1}\frac{2}{j+1}\gamma_{i,j}\right)u^{2p_{0}+1}+\sum_{l=0}^{\left[\frac{n-1}{2}\right]}\sum\limits_{k=0}^{l}\chi_{l,k}\zeta_{2l-2k,2k+1}(u^{2}+u^{2m})^{l+1},
\end{aligned}
\eqno(3.24)$$
where $p_{0}=p+\frac{1}{2}(m-1)(k+1)$. Using (3.17), we have $\omega_{l,0}^{0}=\frac{2}{2l+1}$ and $\chi_{l,0}\neq 0$. This ends the proof. $\diamondsuit$

\vskip 0.2 true cm
\noindent{\bf Lemma 3.5.} For $h\in(0,+\infty)$ and $(m,n)\in D_{1}$, the generating functions of $M(u)$ are the following $\left[\frac{n-1}{2}\right]+\left[\frac{n}{2}\right]^{2}+3\left[\frac{n}{2}\right]+3$ linearly independent functions:
  $$\begin{aligned}
  &u,~u^{3},...,u^{2\left[\frac{n}{2}\right]+1},~u^{m},~u^{m+2},~u^{2m+1},...,u^{2\left[\frac{n}{2}\right]+m},\\
  &...,u^{2\left[\frac{n}{2}\right]m+m},~u^{2}+u^{2m},...,(u^{2}+u^{2m})^{\left[\frac{n-1}{2}\right]+1}.
  \end{aligned}$$

\vskip 0.2 true cm
\noindent{\bf Proof.} The proof for the functions
$$u,~u^3,...,~u^{2t+1},~u^{2}+u^{2m},~(u^{2}+u^{2m})^{2},...,(u^{2}+u^{2m})^{n},~t,n\in\mathbb{N}^{*},$$
are linearly independent is equivalent to showing that
$$1,~u(1+u^{2m-2}),~u^{2},~u^{3}(1+u^{2m-2})^{2},...,~u^{2n-1}(1+u^{2m-2})^{n},~t,n\in\mathbb{N}^{*},$$
are linearly independent. For $t=1$, $n=2$, let
$$\alpha_{0}+\alpha_{1}u(1+u^{2m-2})+\alpha_{2}u^{2}+\alpha_{3}u^{3}(1+u^{2m-2})^{2}=0, \eqno(*_{1})$$
we have $\alpha_{0}=0$ when $u\rightarrow 0$. $(*_{1})$ is equivalent
$$\alpha_{1}(1+u^{2m-2})+\alpha_{2}u+\alpha_{3}u^{2}(1+u^{2m-2})^{2}=0, \eqno(*_{2})$$
we have $\alpha_{1}=0$ when $u\rightarrow 0$. $(*_{2})$ is equivalent
$$\alpha_{2}+\alpha_{3}u(1+u^{2m-2})^{2}=0, \eqno(*_{3})$$
we have $\alpha_{2}=0$ when $u\rightarrow 0$. $(*_{3})$ is equivalent
$$\alpha_{3}(1+u^{2m-2})^{2}=0, \eqno(*_{4})$$
we have $\alpha_{2}=0$ when $u\rightarrow 0$. For $t\geq 2$, $n\geq 3$, we can finish the proof using induction on $n$.

Next, we only need to show that the coeffcients $A_{p}$, $B_{p,k}$ and $K_{l}$ are independent. By Lemma 3.4, we have
  $$
  {\bf G_{1}}=:\frac{\partial{(A_{0},B_{1},...,A_{\left[\frac{n}{2}\right]},B_{0,0},B_{1,0},B_{1,1},...,
  B_{\left[\frac{n}{2}\right],2\left[\frac{n}{2}\right]},K_{0},K_{1},...,K_{\left[\frac{n-1}{2}\right]})}}
  {\partial{(\rho_{0,0},\rho_{2,0},...,\rho_{2\left[\frac{n}{2}\right],0},\gamma_{0,0},\gamma_{2,0},\gamma_{1,1},...,
  \gamma_{0,2\left[\frac{n}{2}\right]},\zeta_{0,1},\zeta_{2,1},...,\zeta_{2\left[\frac{n-1}{2}\right],1})}}
  $$

  $$
    =\begin{pmatrix}
    \begin{smallmatrix}
       &2 &0 &\dots &0 &0 &0 &0&\dots &0 &0 &0 &\dots&0\\
          &0 &\frac{2}{3} &\dots &0&0 &0&0 &\dots &0&0&0&\dots&0\\
          &\dots &\dots &\dots &\dots &\dots &\dots &\dots &\dots &\dots &\dots &\dots &\dots &\dots\\
          &0 &0 &\dots &\frac{2}{2\left[\frac{n}{2}\right]+1} &0 &0 &0 &\dots &0 &0 &0 &\dots&0\\
          &0 &0 &\dots &0 &-2 &0 &0 &\dots &0 &0 &0 &\dots &0\\
          &0 &0 &\dots &0 &0 &-2 &0 &\dots &0 &0 &0 &\dots &0\\
          &0 &0 &\dots &0 &0 &0 &-1 &\dots &0 &0 &0 &\dots &0\\
          &\dots &\dots &\dots &\dots &\dots &\dots &\dots &\dots &\dots &\dots &\dots &\dots &\dots\\
          &0 &0 &\dots &0 &0 &0 &0 &\dots &-\frac{2}{2\left[\frac{n}{2}\right]+1} &0 &0 &\dots &0\\
          &0 &0 &\dots &0 &0 &0 &0 &\dots &0 &\chi_{0,0} &0 &\dots &0\\
          &0 &0 &\dots &0 &0 &0 &0 &\dots &0 &0 &\chi_{1,0} &\dots &0\\
          &\dots &\dots &\dots &\dots &\dots &\dots &\dots &\dots &\dots &\dots &\dots &\dots&\dots\\
          &0 &0 &\dots &0 &0 &0 &0 &\dots &0 &0 &0 &\dots &\chi_{\left[\frac{n-1}{2}\right],0}
       \end{smallmatrix}
       \end{pmatrix}
  $$
  Hence, ${\rm det}{\bf G_{1}}\neq 0$ since $\chi_{k,0}\neq 0$, which implies the independence of the coefficients. This ends the proof. $\diamondsuit$

\vskip 0.2 true cm
\noindent{\bf Proof of the case of} ${\bf (m,n)\in D_{1}}$. By (3.19) and $(u^{2}+u^{2m})^{n}=\sum\limits_{i=0}^{n}\binom{n} i u^{2n-2i+2mi}$, we can get
$$\begin{aligned}
M(u)&=\sum_{p=0}^{\left[\frac{n}{2}\right]}A_{p}u^{2p+1}+\sum_{p=0}^{\left[\frac{n}{2}\right]}\sum_{k=0}^{2p}B_{p,k}u^{2p_{0}+1}
+\sum_{l=0}^{\left[\frac{n-1}{2}\right]}\sum_{i=0}^{l+1}K_{l}\binom{l+1} i u^{2(l+1)-2i+2mi},\\
\end{aligned}
\eqno(3.25)$$
where $p_{0}=p+\frac{1}{2}(m-1)(k+1)$. We claim for $n_{i}\in\mathbb{N},~n_{1}<n_{2}<\dots<n_{k}$, the ordered set $(u^{n_{1}},~u^{n_{2}},~u^{n_{3}},...,u^{n_{k-1}},~u^{n_{k}})$ is an ECT-system on $u\in(0,+\infty)$. We can get the results using induction on $k$. Let $r_{k-1}=n_{k}-n_{1}$. For $k=1$, we have $W[u^{n_{1}}]>0$. For $k=2$, it is equivalent to showing that $1,~u^{r_{1}}$ is an ECT-system, and we can easily get $W\left(1,u^{r_{1}}\right)=r_{1}u^{r_{1}-1}\neq 0$. For $k=3$, it is equivalent to showing that $1,~u^{r_{1}},~u^{r_{2}}$ is an ECT-system. By direct computation,
$$\begin{aligned}
W\left(1,u^{r_{1}},u^{r_{2}}\right)&=
\left|\begin{matrix}
          &1 &u^{r_{1}} &u^{r_{2}}\\
          &0 &r_{1}u^{r_{1}-1} &r_{2}u^{r_{2}-1}\\
          &0 &r_{1}(r_{1}-1)u^{r_{1}-2} &r_{2}(r_{2}-1)u^{r_{2}-2}
          \end{matrix}\right| \\
&=r_{1}r_{2}\left|\begin{matrix}
          &u^{r_{1}-1} &u^{r_{2}}\\
          &(r_{1}-1)u^{r_{1}-2} &(r_{2}-1)u^{r_{2}-2}
          \end{matrix}\right|\\
&=r_{1}r_{2}(r_{2}-r_{1})u^{r_{1}+r_{2}-3}.
\end{aligned}$$
Hence $W\left(1,u^{r_{1}},u^{r_{2}}\right)\neq 0$. Suppose that the results holds for $k\leq p$, then for $k=p+1$, it is equivalent to showing that $1,~u^{r_{1}},~u^{r_{2}},~...,~u^{r_{p}}$ is an ECT-system. By the above computation and assumption, we can get
$$\begin{aligned}
&W\left(1,u^{r_{1}},u^{r_{2}},...,u^{r_{p}}\right)\\
=&
\prod\limits_{i=1}^{p-1}r_{i}\left|\begin{matrix}
          &u^{r_{1}-1} &u^{r_{2}-1} &\dots &u^{r_{p}-1}\\
          &(r_{1}-1)u^{r_{1}-2} &(r_{2}-1)u^{r_{2}-2} &\dots &(r_{p}-1)u^{r_{p}-2}\\
          &\dots &\dots &\dots &\dots\\
          &\prod\limits_{i=1}^{p-1}(r_{1}-i)u^{r_{1}-p} &\prod\limits_{i=1}^{p-1}(r_{2}-i)u^{r_{2}-p} &\dots &\prod\limits_{i=1}^{p-1}(r_{p}-i)u^{r_{p}-p}
          \end{matrix}\right|\neq 0
\end{aligned}$$
Notice (3.25) and $\#\{M(h)=0,h\in(0,+\infty)\}=\#\{M(u)=0,u\in(0,+\infty)\}$, therefore $M(h)$ has at most $\frac{1}{2}\left[\frac{n-1}{2}\right]^{2}+\frac{5}{2}\left[\frac{n-1}{2}\right]+\left[\frac{n}{2}\right]^{2}+3\left[\frac{n}{2}\right]+3$ zeros in $h\in(0,+\infty)$ for $(m,n)\in D_{1}$ using Definition 2.3 and Lemma 2.4.

By Lemma 3.5, $M(u)$ is a linear combination of $\left[\frac{n-1}{2}\right]+\left[\frac{n}{2}\right]^{2}+3\left[\frac{n}{2}\right]+3$ independent functions with arbitrary coefficients for $n<m-1$. All these functions are analytic in $u\in(0,+\infty)$ and strictly positive in this interval. Hence, according to Lemma 2.2 and $\#\{M(h)=0,h\in(0,+\infty)\}=\#\{M(u)=0,u\in(0,+\infty)\}$, there exist coefficients such that $M(h)$ has at least $\left[\frac{n-1}{2}\right]+\left[\frac{n}{2}\right]^{2}+3\left[\frac{n}{2}\right]+2$ zeros in $h\in(0,+\infty)$ for $(m,n)\in D_{1}$. This ends the proof. $\diamondsuit$

\vskip 0.3 true cm
\noindent{\bf Proof of the case of} ${\bf (m,n)\in D_{2} ~and ~(m,n)\in D_{3}}$.
\vskip 0.2 true cm
For these cases, the results of Lemma 3.2 and Lemma 3.3 are also true. Therefore, we have the following results.

\vskip 0.2 true cm
\noindent{\bf Lemma 3.6.} Let $h=u^2+u^{2m}$, then number of zeros of the
 first order Melnikov function $M(h)$ is equal with $M(u)$ in $u\in(0,+\infty)$. For $(m,n)\in D_{2}\cup D_{3}$, $M(u)$ can be written as
  $$M(u)=\sum_{p=0}^{p_{1}}\bar{A}_{p}u^{2p+1}+\sum_{s=0}^{\frac{1}{2}(m-5)}\sum_{k=0}^{2s+1}\bar{B}_{s,k}u^{2p_{2}+1}+\sum_{l=0}^{\left[\frac{n-1}{2}\right]}\bar{K}_{l}(u^2+u^{2m})^{l+1},\eqno(3.26)$$
  where
  $$\begin{aligned}
  p_{1}&=\left[n/2\right]m-(m^2-5m+4)/2,\\
  p_{2}&=\left[n/2\right]m+s-(m^2-5m-k(m-1))/2,\\
  \bar{A}_{p}&=\sum\limits_{\substack{l+j(m-1)=p\\l\leq \left[\frac{n}{2}\right], j\leq l}}\sum\limits_{d=j}^{l}\omega_{l,d}^{j}\rho_{2l-2d,2d}+\sum\limits_{\substack{l+(s+\frac{3}{2})(m-1)=p\\l\leq \left[\frac{n}{2}\right], s\leq l-1}}\sum\limits_{d=0}^{l-1}\vartheta_{l,d}^{s}\rho_{2l-2d-1,2d+1}\\
  &-\sum\limits_{\substack{i+mj+m=2p+1\\i+j\leq 2\left[\frac{n}{2}\right]}}\frac{2}{j+1}\gamma_{i,j},\\ 
  \end{aligned}$$
  $$\begin{aligned}
  \bar{B}_{s,k}&=\sum\limits_{\substack{l+j(m-1)=p_{2}\\l\leq \left[\frac{n}{2}\right], j\leq l}}\sum\limits_{d=j}^{l}\omega_{l,d}^{j}\rho_{2i-2d,2d}+\sum\limits_{\substack{l+(s+\frac{3}{2})(m-1)=p_{2}\\l\leq \left[\frac{n}{2}\right], s\leq l-1}}\sum\limits_{d=0}^{l-1}\vartheta_{l,d}^{s}\rho_{2l-2d-1,2d+1}\\
  &-\sum\limits_{\substack{i+mj+m=2p_{2}+1\\i+j\leq 2\left[\frac{n}{2}\right]}}\frac{2}{j+1}\gamma_{i,j},\\ \bar{K}_{l}&=\sum\limits_{k=0}^{l}\chi_{l,k}\zeta_{2l-2k,2k+1},
  \end{aligned}$$
and $\omega_{l,d}^{j}$, $\vartheta_{l,d}^{s}$ and $\chi_{l,k}$ are constants with $\omega_{l,0}^{0}=\frac{2}{2l+1}$ and $\chi_{l,0}\neq 0$.

\vskip 0.2 true cm
\noindent{\bf Proof.}
Proceeding as the proof of Lemma 3.4, we can obtain
$$\begin{aligned}
M(u)&=\sum_{p=0}^{p_{1}}\left(\sum\limits_{\substack{l+j(m-1)=p\\l\leq \left[\frac{n}{2}\right], j\leq l}}\sum\limits_{d=j}^{l}\omega_{l,d}^{j}\rho_{2l-2d,2d}+\sum\limits_{\substack{l+(s+\frac{3}{2})(m-1)=p\\l\leq \left[\frac{n}{2}\right], s\leq l-1}}\sum\limits_{d=0}^{l-1}\vartheta_{l,d}^{s}\rho_{2l-2d-1,2d+1}\right.\\
  &\left.-\sum\limits_{\substack{i+mj+m=2p+1\\i+j\leq 2\left[\frac{n}{2}\right]}}\frac{2}{j+1}\gamma_{i,j}\right)u^{2p+1}+\sum_{s=0}^{\frac{1}{2}(m-5)}\sum_{k=0}^{2s+1}\left(\sum\limits_{\substack{l+j(m-1)=p_{2}\\l\leq \left[\frac{n}{2}\right], j\leq l}}\sum\limits_{d=j}^{l}\omega_{l,d}^{j}\rho_{2i-2d,2d}\right.\\
  &\left.+\sum\limits_{\substack{l+(s+\frac{3}{2})(m-1)=p_{2}\\l\leq \left[\frac{n}{2}\right], s\leq l-1}}\sum\limits_{d=0}^{l-1}\vartheta_{l,d}^{s}\rho_{2l-2d-1,2d+1}-\sum\limits_{\substack{i+mj+m=2p_{2}+1\\i+j\leq 2\left[\frac{n}{2}\right]}}\frac{2}{j+1}\gamma_{i,j}\right)u^{2p_{2}+1}\\
\end{aligned}\eqno(3.27)$$
$$\begin{aligned}
  &+\sum_{l=0}^{\left[\frac{n-1}{2}\right]}\sum\limits_{k=0}^{l}\mu_{l,k}\zeta_{2l-2k,2k+1}(u^{2}+u^{2m})^{l+1},\quad\quad\quad\quad\quad\quad\quad\quad\quad\quad\quad\quad\quad\quad\quad
\end{aligned}
$$
where $p_{1}=\left[n/2\right]m-(m^2-5m+4)/2$, $p_{2}=\left[n/2\right]m+s-(m^2-5m-k(m-1))/2$.
Using (3.17), we have $\omega_{l,0}^{0}=\frac{2}{2l+1}$ and $\chi_{l,0}\neq 0$. This ends the proof. $\diamondsuit$

\vskip 0.2 true cm
\noindent{\bf Lemma 3.7.} For $h\in(0,+\infty)$, $(m,n)\in D_{2}\cup D_{3}$, the generating functions of $M(u)$ are the following $\left[\frac{n-1}{2}\right]+\left[\frac{n}{2}\right]m-\frac{1}{4}m^{2}+\frac{3}{2}m+\frac{3}{4}$ linearly independent functions:
  $$\begin{aligned}
  &u,~u^{3},~...,~u^{2z+1},~u^{2q+1},~u^{2q+m},~u^{2q+3},u^{2q+m+2},u^{2q+2m+1},\\
  &...,u^{2q+m^{2}-4m},u^{2}+u^{2m},...,(u^{2}+u^{2m})^{\left[\frac{n-1}{2}\right]+1}.
  \end{aligned}$$
  where $z=\left[\frac{n}{2}\right]m-\frac{m^{2}-5m+4}{2}$, $q=\left[\frac{n}{2}\right]m-\frac{m^{2}-5m}{2}$.

\vskip 0.2 true cm
\noindent{\bf Proof.}
  We only need to prove that the coefficients $\bar{A}_{p}$, $\bar{B}_{s,k}$ and $\bar{K}_{l}$ are independent. From Lemma 3.6, we can get
  \begin{scriptsize}$$
  {\bf G_{2}}=:\frac{\partial{(\bar{A}_{0},\bar{A}_{1},...,\bar{A}_{\frac{m-3}{2}},\bar{A}_{\frac{m-1}{2}},...,
  \bar{A}_{n_{1}},\bar{A}_{n_{2}},...,\bar{A}_{n_{3}},
  \bar{B}_{0,0},\bar{B}_{0,1},...,\bar{B}_{\frac{m-5}{2},0},...,\bar{B}_{\frac{m-5}{2},m-4},\bar{K}_{0},\bar{K}_{1},...,\bar{K}_{\left[\frac{n-1}{2}\right]})}}
  {\partial{(\rho_{0,0},\rho_{2,0},...,\rho_{m-3,0},\gamma_{0,0},...,
  \gamma_{k_{1},0},\gamma_{k_{2},1},...,\gamma_{m-3,k_{3}},\gamma_{1,k_{4}},
  \gamma_{0,k_{5}},...,\gamma_{m-4,k_{6}},...,\gamma_{0,k_{7}},\zeta_{0,1},\zeta_{2,1},...,\zeta_{2\left[\frac{n-1}{2}\right],1})}}
  $$\end{scriptsize}
  $$\tiny{
  \addtocounter{MaxMatrixCols}{15}
    =\begin{pmatrix}
    \begin{smallmatrix}
          &2 &0 &\dots &0 &0 &\dots &0 &0 &\dots &0 &0 &0 &\dots &0 &\dots &0 &0 &0 &\dots &0\\
          &0 &\frac{2}{3} &\dots &0 &0 &\dots &0 &0 &\dots &0 &0 &0 &\dots &0 &\dots &0 &0 &0 &\dots &0\\
          &\dots &\dots &\dots &\dots &\dots &\dots &\dots &\dots &\dots &\dots &\dots &\dots &\dots &\dots &\dots &\dots &\dots &\dots &\dots &\dots  \\
          &0 &0 &\dots &\frac{2}{m-2} &0 &\dots &0 &0 &\dots &0 &0 &0 &\dots &0 &\dots &0 &0 &0 &\dots &0\\
          &0 &0 &\dots &0 &-2 &\dots &0 &0 &\dots &0 &0 &0 &\dots &0 &\dots &0 &0 &0 &\dots &0\\
          &\dots &\dots &\dots &\dots &\dots &\dots &\dots &\dots &\dots &\dots &\dots &\dots &\dots &\dots &\dots &\dots &\dots &\dots &\dots &\dots  \\
          &0 &0 &\dots &0 &0 &\dots &-2 &0 &\dots &0 &0 &0 &\dots &0 &\dots &0 &0 &0 &\dots &0\\
          &0 &0 &\dots &0 &0 &\dots &0 &-1 &\dots &0 &0 &0 &\dots &0 &\dots &0 &0 &0 &\dots &0\\
          &\dots &\dots &\dots &\dots &\dots &\dots &\dots &\dots &\dots &\dots &\dots &\dots &\dots &\dots &\dots &\dots &\dots &\dots &\dots &\dots  \\
          &0 &0 &\dots &0 &0 &\dots &0 &0 &\dots &-\frac{2}{2\left[\frac{n}{2}\right]-m+4} &0 &0 &\dots &0 &\dots &0 &0 &0 &\dots &0\\
          &0 &0 &\dots &0 &0 &\dots &0 &0 &\dots &0 &-\frac{2}{2\left[\frac{n}{2}\right]-m+5} &0 &\dots &0 &\dots &0 &0 &0 &\dots &0\\
          &0 &0 &\dots &0 &0 &\dots &0 &0 &\dots &0 &0 &-\frac{2}{2\left[\frac{n}{2}\right]-m+6} &\dots &0 &\dots &0 &0 &0 &\dots &0\\
          &\dots &\dots &\dots &\dots &\dots &\dots &\dots &\dots &\dots &\dots &\dots &\dots &\dots &\dots &\dots &\dots &\dots &\dots &\dots &\dots  \\
          &0 &0 &\dots &0 &0 &\dots &0 &0 &\dots &0 &0 &0 &\dots &-\frac{2}{2\left[\frac{n}{2}\right]-m+5} &\dots &0 &0 &0 &\dots &0\\
          &\dots &\dots &\dots &\dots &\dots &\dots &\dots &\dots &\dots &\dots &\dots &\dots &\dots &\dots &\dots &\dots &\dots &\dots &\dots &\dots  \\
          &0 &0 &\dots &0 &0 &\dots &0 &0 &\dots &0 &0 &0 &\dots &0 &\dots &-\frac{2}{2\left[\frac{n}{2}\right]+1} &0 &0 &\dots &0\\
          &0 &0 &\dots &0 &0 &\dots &0 &0 &\dots &0 &0 &0 &\dots &0 &\dots &0 &\chi_{0,0} &0 &\dots &0\\
          &0 &0 &\dots &0 &0 &\dots &0 &0 &\dots &0 &0 &0 &\dots &0 &\dots &0 &0 &\chi_{1,0} &\dots &0\\
          &\dots &\dots &\dots &\dots &\dots &\dots &\dots &\dots &\dots &\dots &\dots &\dots &\dots &\dots &\dots &\dots &\dots &\dots &\dots &\dots  \\
          &0 &0 &\dots &0 &0 &\dots &0 &0 &\dots &0 &0 &0 &\dots &0 &\dots &0 &0 &0 &\dots &\chi_{\left[\frac{n-1}{2}\right],0}
       \end{smallmatrix}
       \end{pmatrix}
  }$$
where
$$\begin{aligned}
 &n_{1}=\left[n/2\right]+\frac{m-1}{2}, ~~&n_{2}=\left[n/2\right]+\frac{m+1}{2}, ~~&n_{3}=\left[n/2\right]m-\frac{m^{2}-5m+4}{2},\\
 &k_{1}=2\left[n/2\right], ~~&k_{2}=2\left[n/2\right]-m+2, ~~&k_{3}=2\left[n/2\right]-m+3,\\
 &k_{4}=2\left[n/2\right]-m+4, ~~&k_{5}=2\left[n/2\right]-m+5, ~~&k_{6}=2\left[n/2\right]-m+4,\\
 &k_{7}=2\left[n/2\right].
\end{aligned}$$
So ${\rm det}{\bf G_{2}}\neq0$ since $\chi_{l,0}\neq 0$, which implies the independence of the coefficients. This ends the proof. $\diamondsuit$

Similar to the prove for the case of $(m,n)\in D_{1}$, we can get the results of Theorem 1.1 (ii).

\vskip 0.5 true cm
\centerline{\bf{ $\S$4}. Proof of Theorem 1.2.}
\vskip 0.5 true cm

In this section, we will consider the case for $m=2k$.
Denote $\rho_{i,j}=b_{i,j}^{+}+\frac{i+1}{j}a_{i+1,j-1}^{+}-b_{i,j}^{-}-\frac{i+1}{j}a_{i+1,j-1}^{-}$, $\gamma_{i,j}=a_{i,j}^{+}-a_{i,j}^{-}$ and $\zeta_{i,j}=b_{i,j}^{-}+\frac{i+1}{j}a_{i+1,j-1}^{-}$, and it is easy to show that the coefficients $\rho_{i,j}$, $\gamma_{i,j}$ and $\zeta_{i,j}$ are independent. In fact, let
$$\overline{{\bf M}}:=\frac{\partial{(\rho_{0,0},\rho_{1,0},\rho_{0,1},...,\gamma_{0,0},\gamma_{1,0},\gamma_{0,1},...\zeta_{0,0},\zeta_{1,0},\zeta_{0,1}...,\zeta_{0,n})}}
{\partial{(b_{0,0}^{+},b_{1,0}^{+},b_{0,1}^{+},...,a_{0,0}^{+},a_{1,0}^{+},a_{0,1}^{+},...,,b_{0,0}^{-},b_{1,0}^{-},b_{0,1}^{-}...,b_{0,n}^{-})}},$$
according to the expression of $\rho_{i,j}$, $\gamma_{i,j}$ and $\zeta_{i,j}$ we have ${\rm deg}~\overline{{\bf M}}=1$.

By direct computation, we have $J_{1,0}=I_{1,0}=J_{1,1}=I_{1,1}=0$,
$$J_{0,0}=-I_{0,0}=2u,~~J_{0,1}=2\int^{u}_{0}\sqrt{h-x^{2}}dx,~~I_{0,1}=\pi h-2\int^{u}_{0}\sqrt{h-x^{2}}dx.$$
Thus, similar to the proof of Lemma 3.2 and Lemma 3.3, we can obtain the following Lemma.

\vskip 0.2 true cm
\noindent{\bf Lemma 4.1.} Let $h=u^{2}+u^{2m}$, then number of zeros of the first order Melnikov function $M(h)$ is equal with $M(u)$ in $u\in(0,+\infty)$. And $M(u)$ can be writen as
$$\begin{aligned}
M(u)&=\sum_{l=0}^{\left[\frac{n}{2}\right]}\sum_{j=0}^{l}C_{l,j}u^{2l+2j(m-1)+1}+\sum_{l=0}^{\left[\frac{n-1}{2}\right]}\sum_{j=0}^{l}D_{l,j}u^{2l+2j(m-1)+m+1}\\
&+\sum_{l=0}^{\left[\frac{n-1}{2}\right]}E_{l}(u^{2}+u^{2m})^{l+1}+\sum_{l=0}^{\left[\frac{n-1}{2}\right]}F_{l}(u^{2}+u^{2m})^{l+1}\int_{0}^{1/\sqrt{1+u^{2m-2}}}\sqrt{1-t^{2}}dt,
\end{aligned}\eqno(4.1)$$
where
$$\begin{aligned}
C_{l,0}&=\sum_{d=0}^{l}\omega_{l,d}^{0}\rho_{2l-2d,2d},~~D_{l,0}=-2\gamma_{2l+1,0},\\
C_{l,j}&=\sum_{d=j}^{l}\omega_{l,d}^{j}\rho_{2l-2d,2d}-\frac{1}{j}\gamma_{2l-2j+1,2j-1},~j\geq 1,\\
D_{l,j}&=\sum_{d=0}^{l}\nu_{l,d}^{j-1}\rho_{2l-2d,2d+1}-\frac{2}{2j+1}\gamma_{2l-2j+1,2j},~j\geq 1,\\
E_{l}&=\sum_{d=0}^{l}\pi\nu_{l,d}^{0}\zeta_{2l-2d,2d+1},~~F_{l}=\sum_{d=0}^{l}2\nu_{l,d}^{0}\rho_{2l-2d,2d+1},
\end{aligned}$$
$\omega_{l,d}^{j}$ and $\nu_{l,d}^{j}$ are constants with $\omega_{l,0}^{0}=\frac{2}{2l+1}$ and $\nu_{l,0}^{0}\neq 0$.

The Taylor expansion of function $\int_{0}^{1/\sqrt{1+u^{2m-2}}}\sqrt{1-t^{2}}dt$ in the variable $u$, around $u=0$, is
$$\begin{aligned}
\int_{0}^{1/\sqrt{1+u^{2m-2}}}\sqrt{1-t^{2}}dt=&\frac{\pi}{4}-\frac{1}{3}u^{3m-3}+\frac{2}{5}u^{5m-5}-\frac{3}{7}u^{7m-7}+\frac{4}{9}u^{9m-9}+\dots\\
=&\frac{\pi}{4}+\sum_{k=1}^{\infty}(-1)^{k}\frac{k}{2k+1}u^{(2k+1)(m-1)},
\end{aligned}$$
thus we have
$$(u^{2}+u^{2m})^{l}\int_{0}^{1/\sqrt{1+u^{2m-2}}}\sqrt{1-t^{2}}dt=\frac{\pi}{4}(u^{2}+u^{2m})^{l}+\sum_{i=0}^{l}\sum_{k=1}^{\infty}\binom{l} i (-1)^{k}\frac{k}{2k+1}u^{2l+(2i+2k+1)(m-1)}.$$
Combining Lemma 4.1, for $n\leq 2m$, the function of $M(u)$ is analytic at $u=0$ in $u$ and has an expansion as follows:
$$M(u)=\alpha_{1}g_{1}(u)+\alpha_{2}g_{2}(u)+...+\alpha_{i}g_{i}(u)+\alpha_{i+1}g_{i+1}(u)+...,~g_{i+1}(u)=o\left(g_{i}(u)\right),\eqno(4.2)$$
which is convergent for small $u$ $(0<u\ll 1)$.
\vskip 0.2 true cm
\noindent{\bf Lemma 4.2.} The function of $M(u)$ in (4.2) can be expressed as
(i)~~when $n<m$, we have
  $$\begin{aligned}
  M(u)&=\Phi_{1}(u)+\Psi_{1}(u)+\Phi_{2}(u)+\Upsilon_{1}(u)+\Psi_{2}(u)+\Phi_{3}(u)+...+\Phi_{\frac{n}{2}}(u)\\
  &+\Upsilon_{\frac{n}{2}-1}(u)+\Psi_{\frac{n}{2}}(u)+\Phi_{\frac{n}{2}+1}(u)+\Upsilon_{\frac{n}{2}}(u)+o(u^{(n+1)m-1}),~n~ is~ even,
\end{aligned}\eqno(4.3)$$
  $$\begin{aligned}
  M(u)&=\Phi_{1}(u)+\Psi_{1}(u)+\Phi_{2}(u)+\Upsilon_{1}(u)+\Psi_{2}(u)+\Phi_{3}(u)+...+\Psi_{\frac{n-1}{2}}(u)\\
  &+\Phi_{\frac{n+1}{2}}(u)+\Upsilon_{\frac{n-1}{2}}(u)+\Psi_{\frac{n+1}{2}}(u)+\Upsilon_{\frac{n+1}{2}}(u)+o(u^{(n+2)m-1}),~n~is~odd,
\end{aligned}\eqno(4.4)$$

  where
  $$\begin{aligned}
  \xi_{2l+1}^{(1)}&=C_{l,0},~~\xi_{2l+2}^{(1)}=E_{l}+\frac{\pi}{4}F_{l},~\xi_{l-k+1}^{(k)}=C_{l,k-1}~(k\geq 2),\\
  \tau_{l}^{(1)}&=-2\gamma_{2l+1,0},~\chi_{k}=F_{k-1}~(k\geq 1),\\
  \tau_{l-k+1}^{(k)}&=D_{l,k-1}+ \left(\sum_{i=0}^{min\{l+1,k-2\}}\binom{l+1} i (-1)^{k-1-i}\frac{k-1-i}{2(k-1-i)+1}\right)F_{l},~k\geq 2,\\
  \Phi_{1}(u)&=\sum_{l=0}^{\left[\frac{n}{2}\right]}\xi_{2l+1}^{(1)}u^{2l+1}+\sum_{l=0}^{\left[\frac{n-1}{2}\right]}\xi_{2l+2}^{(1)}(u^{2}+u^{2m})^{l+1},\\
  \Phi_{k}(u)&=\sum_{l=k-1}^{\left[\frac{n}{2}\right]}\xi_{l-k+1}^{(k)}u^{2l+2(k-1)(m-1)+1},~k=2,3,...,\left[n/2\right]+1,\\
  \Psi_{k}(u)&=\sum_{l=k-1}^{\left[\frac{n-1}{2}\right]}\tau_{l-k+1}^{(k)}u^{2l+2(k-1)(m-1)+m+1},~k=1,2,...,\left[(n-1)/2\right]+1,\\
  \Upsilon_{k}(u)&=\chi_{k}u^{(2k+1)m-1}\left(\sum_{i=0}^{k-1}\binom {k} i (-1)^{k-i}\frac{k-i}{2(k-i)+1}\right.\\
  & \left.+\left(\sum_{j=1}^{\left[\frac{n+1}{2}\right]-k}\sum_{i=0}^{k}\binom {k} i (-1)^{k+j-i}\frac{k+j-i}{2(k+j-i)+1}\right)u^{2j(m-1)}\right),\\
  &k=1,2,...,\left[(n-1)/2\right]+1,\\
  \end{aligned}$$
  $C_{l,i}$, $D_{l,j}$, $E_{l}$ and $F_{l}$ are shown in (4.1).
(ii)~~ When $m\leq n<2m-2$, we have
  $$\begin{aligned}
  M(u)&=\bar{\Phi}_{1}(u)+\bar{\Psi}_{1}(u)+\bar{\Phi}_{2}(u)+...+\bar{\Phi}_{\frac{n-m}{2}+3}(u)+\bar{\Upsilon}_{1}(u)+\bar{\Psi}_{\frac{n-m}{2}+3}(u)\\
  &+\bar{\Phi}_{\frac{n-m}{2}+4}(u)+\bar{\Upsilon}_{2}(u)+\bar{\Psi}_{\frac{n-m}{2}+4}(u)...+\bar{\Psi}_{\frac{m}{2}+1}(u)+\bar{\Lambda}_{1}(u)\\
  &+\bar{\Phi}_{\frac{m}{2}+2}(u)+\bar{\Upsilon}_{\frac{2m-n}{2}}(u)+...+\bar{\Lambda}_{\frac{n-m}{2}}(u)+\bar{\Phi}_{\frac{n}{2}+1}+\bar{\Upsilon}_{\frac{m}{2}-1}(u)\\
  &+\bar{\Lambda}_{\frac{n-m}{2}+1}(u)+o(u^{(n+2)m-1}),~n~is~even,
\end{aligned}\eqno(4.5)$$
  $$\begin{aligned}
  M(u)&=\bar{\Phi}_{1}(u)+\bar{\Psi}_{1}(u)+\bar{\Phi}_{2}(u)+...+\bar{\Phi}_{\frac{n-m+5}{2}}(u)+\bar{\Upsilon}_{1}(u)+\bar{\Psi}_{\frac{n-m+5}{2}}(u)\\
  &+\bar{\Phi}_{\frac{n-m+7}{2}}(u)+\bar{\Upsilon}_{2}(u)+\bar{\Psi}_{\frac{n-m+7}{2}}(u)...+\bar{\Psi}_{\frac{m}{2}+1}(u)+\bar{\Lambda}_{1}(u)\\
  &+\bar{\Phi}_{\frac{m}{2}+2}(u)+\bar{\Upsilon}_{\frac{2m-n-3}{2}}(u)+...+\bar{\Upsilon}_{\frac{m}{2}-1}(u)+\bar{\Psi}_{\frac{n+1}{2}}(u)+\bar{\Lambda}_{\frac{n-m+1}{2}}(u)\\
  &+\bar{\Upsilon}_{\frac{m}{2}}(u)+o(u^{(n+2)m-1}),~n~is~odd,
\end{aligned}\eqno(4.6)$$
  where
  $$
  \begin{aligned}
  \bar{\Phi}_{1}(u)&=\sum_{l=0}^{\left[\frac{n}{2}\right]}\bar{\xi}_{2l+1}^{(1)}u^{2l+1}+\sum_{l=0}^{\left[\frac{n-1}{2}\right]}\bar{\xi}_{2l+2}^{(1)}(u^{2}+u^{2m})^{l+1},\\
  \bar{\Phi}_{k}(u)&=\begin{cases}\sum\limits_{l=t_{1}}^{\left[\frac{n}{2}\right]}\bar{\xi}_{l-t_{1}}^{(k)}u^{2l+2(k-1)(m-1)+1},~k=2,3,...,t_{1}+1,~t_{1}=\left[(n-1)/2\right]-m/2+2,\\
                  \sum\limits_{l=k-1}^{\left[\frac{n}{2}\right]}\bar{\xi}_{l-k+1}^{(k)}u^{2l+2(k-1)(m-1)+1},~k=t_{1}+2,...,\left[n/2\right]+1,\\
                  \end{cases}\\
  \bar{\Psi}_{k}(u)&=\begin{cases}\sum\limits_{l=t_{2}}^{\left[\frac{n-1}{2}\right]}\bar{\tau}_{l-t_{2}}^{(k)}u^{2l+2(k-1)(m-1)+m+1},~k=1,2,...,t_{2}+1,~t_{2}=\left[n/2\right]-m/2+1,\\
                  \sum\limits_{l=k-1}^{\left[\frac{n-1}{2}\right]}\bar{\tau}_{l-k+1}^{(k)}u^{2l+2(k-1)(m-1)+m+1},~k=t_{2}+2,...,\left[(n-1)/2\right]+1,\\
                  \end{cases}\\
  \bar{\Upsilon}_{k}(u)&=\bar{\chi}_{k}u^{2t_{2}(m-1)+n+2km-\delta_{n}} \left(\sum_{i=0}^{t_{2}-1+k}\binom {t_{2}+k} i (-1)^{t_{2}+k-i}\frac{t_{2}+k-i}{2(t_{2}+k-i)+1}\right.\\
  & \left.+\sum_{j=1}^{\frac{m}{2}+\delta_{n}-k}\sum_{i=0}^{t_{2}+k}\binom {t_{2}+k} i (-1)^{t_{2}+k+j-i}\frac{t_{2}+k+j-i}{2(t_{2}+k+j-i)+1}u^{2j(m-1)}\right),\\
  &\quad k=1,2,...,m/2+\delta_{n},\\
  \bar{\Lambda}_{k}(u)&=\bar{\eta}_{k}u^{(2k+m)m-1} \left(\sum_{i=0}^{k}\binom {k} i (-1)^{\frac{m}{2}+k-i}\frac{m/2+k-i}{2(m/2+k-i)+1}\right.\\
  & \left.+\sum_{j=1}^{\left[\frac{n}{2}\right]-\frac{m}{2}+1-k}\sum_{i=0}^{k}\binom {k} i (-1)^{\frac{m}{2}+k+j-i}\frac{m/2+k+j-i}{2(m/2+k+j-i)+1}u^{2j(m-1)}\right),\\
  &\quad k=1,2,...,\left[n/2\right]-m/2+1,\\
  \end{aligned}$$
  $$\begin{aligned}
  \bar{\xi}_{2l+1}^{(1)}&=C_{l,0}-2\gamma_{2l-m+1,0},~\bar{\xi}_{2l+2}^{(1)}=E_{l}+\frac{\pi}{4}F_{l},\\
  \bar{\tau}_{l-t_{2}}^{(1)}&=-2\gamma_{2l+1,0}+C_{l-\frac{m}{2}+1,1},\\
  \bar{\xi}_{l-t_{1}}^{(k)}&= \left(\sum_{i=0}^{min\{l-\frac{m}{2}+1,k-2\}}\binom{l-\frac{m}{2}+1} i (-1)^{k-1-i}\frac{k-1-i}{2(k-1-i)+1}\right)F_{l-\frac{m}{2}}\\
  &+C_{l,k-1}+D_{l-\frac{m}{2},k-1},~k=2,3,...,t_{1}+1,\\
  \bar{\tau}_{l-t_{2}}^{(k)}&=D_{l,k-1}+C_{l-\frac{m}{2}+1,k}+ \left(\sum_{i=0}^{min\{l+1,k-2\}}\binom{l+1} i (-1)^{k-1-i}\frac{k-1-i}{2(k-1-i)+1}\right)F_{l},\\
  &\quad k=2,3,...,t_{2}+1,\\
  \bar{\chi}_{k}&=F_{t_{3}},~t_{3}=\frac{1}{2}(n-m-\delta_{n}-1+2k),~k\geq 1,\\
  \bar{\eta}_{k}&=F_{k-1},~k\geq 1
  \end{aligned}$$
  the expression of $\bar{\xi}_{l-k+1}^{(k)}$ and $\bar{\tau}_{l-k+1}^{(k)}$ are same as $\bar{\xi}_{l-t_{1}}^{(k)}$ and $\bar{\tau}_{l-t_{2}}^{(k)}$ respectively.
(iii)~~When $n=2m-2$, $2m-1$ and $2m$, we have
  $$\begin{aligned}
  M(u)&=\tilde{\Phi}_{1}(u)+\tilde{\Psi}_{1}(u)+\tilde{\Phi}_{2}(u)+...+\tilde{\Psi}_{\frac{n-m}{2}+2}(u)+\tilde{\Lambda}_{1}(u)+\tilde{\Phi}_{\frac{n-m}{2}+3}(u)\\
  &+\tilde{\Upsilon}_{1}(u)+\tilde{\Psi}_{\frac{n-m}{2}+3}(u)+...+\tilde{\Psi}_{\frac{n}{2}}(u)+\tilde{\Lambda}_{\frac{m}{2}-1}(u)+\tilde{\Phi}_{\frac{n}{2}+1}(u)\\
  &+\tilde{\Upsilon}_{\frac{m}{2}-1}(u)+\tilde{\Gamma}_{\frac{m}{2}}(u)+o(u^{(n+2)m-1}),~n~is~even,
\end{aligned}\eqno(4.7)$$
  $$\begin{aligned}
  M(u)&=\tilde{\Phi}_{1}(u)+\tilde{\Psi}_{1}(u)+\tilde{\Phi}_{2}(u)+...+\tilde{\Phi}_{\frac{n-m+5}{2}}(u)+\tilde{\Upsilon}_{1}(u)+\tilde{\Psi}_{\frac{n-m+5}{2}}(u)\\
  &+\tilde{\Lambda}_{1}(u)+\tilde{\Phi}_{\frac{n-m+7}{2}}(u)+...+\tilde{\Phi}_{\frac{n+1}{2}}(u)+\tilde{\Upsilon}_{\frac{m}{2}-1}(u)+\tilde{\Psi}_{\frac{n+1}{2}}(u)\\
  &+\tilde{\Lambda}_{\frac{m}{2}-1}(u)+\tilde{\Upsilon}_{\frac{m}{2}}(u)+o(u^{(n+2)m-1}),~n~is~odd,
\end{aligned}\eqno(4.8)$$
  where
  $$\begin{aligned}
  \tilde{\Phi}_{1}=&\sum_{l=0}^{\left[\frac{n}{2}\right]}\tilde{\xi}_{2l+1}^{(1)}u^{2l+1}+\sum_{l=0}^{\left[\frac{n-1}{2}\right]}\tilde{\xi}_{2l+2}^{(1)}(u^{2}+u^{2m})^{l+1},\\
  \tilde{\Phi}_{k}=&\begin{cases}\sum\limits_{l=t_{1}}^{\left[\frac{n}{2}\right]}\tilde{\xi}_{l-t_{1}}^{(k)}u^{2l+2(k-1)(m-1)+1},~k=2,3,...,t_{1}+1,~t_{1}=\left[(n-1)/2\right]-m/2+2,\\
                   \sum\limits_{l=k-1}^{\left[\frac{n}{2}\right]}\tilde{\xi}_{l-k+1}^{(k)}u^{2l+2(k-1)(m-1)+1},~k=t_{1}+2,...,\left[n/2\right]+1,\\
                   \end{cases}\\
  \tilde{\Psi}_{k}=&\begin{cases}\sum\limits_{l=t_{2}}^{\left[\frac{n-1}{2}\right]}\tilde{\tau}_{l-t_{2}}^{(k)}u^{2l+2(k-1)(m-1)+m+1},~k=1,2,...,t_{2}+1,~t_{2}=\left[n/2\right]-m/2+1,\\
                   \sum\limits_{l=k-1}^{\left[\frac{n-1}{2}\right]}\tilde{\tau}_{l-k+1}^{(k)}u^{2l+2(k-1)(m-1)+m+1},~k=t_{2}+2,...,\left[(n-1)/2\right]+1,\\
                   \end{cases}\\
  \tilde{\Upsilon}_{k}=&\begin{cases}\bar{\chi}_{k}u^{2t_{2}(m-1)+n+2km-\delta_{n}} \left(\sum\limits_{i=0}^{t_{2}-1+k}\binom {t_{2}+k} i (-1)^{t_{2}+k-i}\frac{t_{2}+k-i}{2(t_{2}+k-i)+1}\right.\\
  \left.+\sum\limits_{j=1}^{\frac{m}{2}+\delta_{n}-k}\left(\sum\limits_{i=0}^{t_{2}+k}\binom {t_{2}+k} i (-1)^{t_{2}+k+j-i}\frac{t_{2}+k+j-i}{2(t_{2}+k+j-i)+1}\right)u^{2j(m-1)}\right),\\
  k=1,2,...,\frac{m}{2}-1,~if~n=2m-2;\\
  k=1,2,...,\frac{m}{2}+\delta_{n}-1,~if~n=2m-1,~2m,\\
  \tilde{\chi}_{k}u^{2t_{2}(m-1)+n+2km-\delta_{n}}\left(\sum\limits_{i=0}^{m-1}\binom {m} i (-1)^{m-i}\frac{m-i}{2(m-i)+1}+\sum\limits_{i=0}^{1}\binom{1} i (-1)^{m+i-i}\frac{m+1-i}{2(m+i-i)+1}\right.\\
  \left.+\sum\limits_{j=1}^{\frac{m}{2}+\delta_{n}-k}\left(\sum\limits_{i=0}^{m}\binom {m} i (-1)^{m+j-i}\frac{m+j-i}{2(m+j-i)+1}\right.\right.\\
  \left.\left.+\sum\limits_{i=0}^{1}\binom{1} i (-1)^{m+i+j-i}\frac{m+1+j-i}{2(m+1+j-i)+1}\right)u^{2j(m-1)}\right),\\
  k=\frac{m}{2}+\delta_{n},~if~n=2m-1,n=2m,\\
  \end{cases}\\
   \end{aligned}$$
  $$\begin{aligned}
  \tilde{\Lambda}_{k}&=\tilde{\eta}_{k}u^{2m(t_{1}+k)-2\left[\frac{n-1}{2}\right]+n+\delta_{n}-2} \left(\sum_{i=0}^{t_{1}-\frac{m}{2}+k}\binom {t_{1}-\frac{m}{2}+k} i (-1)^{t_{1}+k-i}\frac{t_{1}+k-i}{2(t_{1}+k-i)+1}\right.\\
  &\left.+\sum_{j=1}^{\frac{m}{2}-\delta_{n}-1-k}\sum_{i=0}^{t_{1}-\frac{m}{2}+k}\binom{t_{1}-\frac{m}{2}+k} i (-1)^{t_{1}+k+j-i}\frac{t_{1}+k+j-i}{2(t_{1}+k+j-i)+1}u^{2j(m-1)}\right),\\
  &\quad k=1,2,...,m/2-1-\delta_{n},\\
  \tilde{\xi}_{2l+1}^{(1)}&=C_{l,0}+C_{l-m+1,1}-2\gamma_{2l-m+1,0},~\tilde{\xi}_{2l+2}^{(1)}=E_{l}+\frac{\pi}{4}F_{l},\\
  \tilde{\xi}_{l-t_{1}}^{(k)}&=\left(\sum_{i=0}^{min\{l-\frac{m}{2}+1,k-2\}}\binom{l-\frac{m}{2}+1} i (-1)^{k-1-i}\frac{k-1-i}{2(k-1-i)+1}\right)F_{l-\frac{m}{2}}\\
  &+C_{l,k-1}+D_{l-\frac{m}{2},k-1},~k=2,3,...,t_{1}+1,\\
  \tilde{\tau}_{l-t_{2}}^{(1)}&=C_{l-\frac{m}{2}+1,1}+D_{l-m+1,1}-2\gamma_{2l+1,0},\\
  \tilde{\tau}_{l-t_{2}}^{(k)}&=\begin{cases}\left(\sum\limits_{i=0}^{min\{l+1,k-2\}}\binom{l+1} i (-1)^{k-1-i}\frac{k-1-i}{2(k-1-i)+1}\right)F_{l}+D_{l,k-1}+C_{l-\frac{m}{2}+1,k},\\
  l\leq \left[(n-1)/2\right],~k\leq t_{2}+1,~if~ n=2m-2;\\
  l\leq \left[(n-1)/2\right]-1,~k\leq t_{2}+1,~if ~n=2m-1,~2m;\\
  \left(\sum\limits_{i=0}^{min\{l+1,k-2\}}\binom{l+1} i (-1)^{k-1-i}\frac{k-1-i}{2(k-1-i)+1}\right)F_{l}+ \left(\sum\limits_{i=0}^{1}(-1)^{k-i}\frac{k-i}{2(k-i)+1}\right)F_{0}\\
  +D_{l,k-1}+C_{l-\frac{m}{2}+1,k},~l=\left[(n-1)/2\right],~if~n=2m-1,~2m,\\
  \end{cases}\\
  \tilde{\chi}_{k}&=F_{t_{3}},~t_{3}=\frac{1}{2}(n-m-\delta_{n}-1+2k),~k\geq 1,\\
  \tilde{\eta}_{k}&=F_{t_{4}},~~t_{4}=n/2-m+1+\delta_{n}+k,~k\geq 1,
  \end{aligned}$$
the expression of $\xi_{l-k+1}^{(k)}$ and $\tau_{l-k+1}^{(k)}$ are same as $\xi_{l-t_{1}}^{(k)}$ and $\tau_{l-t_{2}}^{(k)}$ respectively.

\vskip 0.2 true cm
\noindent{\bf Proof of the case of} ${\bf (m,n)\in D_{4}}.$
\vskip 0.2 true cm

First, we will prove the lower bound of $M(h)$. Without loss of generality, we only prove the case of $n$ is even, the case of $n$ is odd can be shown in a similar way.

We claim the coefficients of $\Phi_{k}$, $\Psi_{k}$, $\Upsilon_{k}$ and $\Gamma_{k}$ are independent. In fact for $n<m$,
$${\bf \overline{G}_{1}}=:\frac{\partial{(\xi_{0}^{(1)},\xi_{1}^{(1)},...,\tau_{0}^{(1)},\tau_{1}^{(1)},...,\chi_{1},\tau_{0}^{(2)},...,\xi_{0}^{(\frac{n}{2})},\xi_{1}^{(\frac{n}{2})},...,\chi_{\frac{n}{2}})}}
{\partial{(\rho_{0,0},\zeta_{0,1},...,\gamma_{1,0},\gamma_{3,0},...,\rho_{0,1},\gamma_{1,2},...,\gamma_{1,n-3},\gamma_{3,n-3},...,\rho_{n-2,1})}}$$

$$
    =\begin{pmatrix}
    \begin{smallmatrix}
          &2 &0 &\dots &0 &0 &\dots &0 &0 &\dots&0 &0 &\dots &0\\
          &0 &\pi\nu_{0,0}^{0} &\dots &0&0 &\dots &\frac{1}{2}\pi\nu_{0,0}^{0} &0 &\dots &0 &0 &\dots &0\\
          &\dots &\dots &\dots &\dots &\dots &\dots &\dots &\dots &\dots &\dots &\dots &\dots &\dots\\
          &0 &0 &\dots &-2 &0 &\dots &0 &0 &\dots &0 &0 &\dots &0\\
          &0 &0 &\dots &0 &-2 &\dots &0 &0 &\dots &0 &0 &\dots &0\\
          &\dots &\dots &\dots &\dots &\dots &\dots &\dots &\dots &\dots &\dots &\dots &\dots &\dots\\
          &0 &0 &\dots &0 &0 &\dots &2\nu_{0,0}^{0} &0 &\dots &0 &0 &\dots &0\\
          &0 &0 &\dots &0 &0 &\dots &0 &-\frac{2}{3} &\dots &0 &0 &\dots &0\\
          &\dots &\dots &\dots &\dots &\dots &\dots &\dots &\dots &\dots &\dots &\dots &\dots &\dots\\
          &0 &0 &\dots &0 &0 &\dots &0 &0 &\dots &\frac{2}{2-n} &0 &\dots &0\\
          &0 &0 &\dots &0 &0 &\dots &0 &0 &\dots &0 &\frac{2}{2-n} &\dots &0\\
          &\dots &\dots &\dots &\dots &\dots &\dots &\dots &\dots &\dots &\dots &\dots &\dots &\dots\\
          &0 &0 &\dots &0 &0 &\dots &0 &0 &\dots &0 &0 &\dots &2\nu_{\frac{n}{2}-1,0}^{0}
       \end{smallmatrix}
       \end{pmatrix}
  $$
hence $\rm{det}{\bf \overline{G}_{1}}\neq 0$ since $\nu_{l,0}^{0}\neq 0$, which implies the independence of the coefficients. Therefore, there exist $\rho_{i,j}$, $\gamma_{i,j}$ and $\zeta_{i,j}$ such that the lower bound of $Z(m,n)$ in Theorem 1.2 (i) can be reached.

Next, we will prove the upper bound of $M(h)$. Denote
$$M_{0}(u)=:(u^{2}+u^{2m})\left(\sum\limits_{l=0}^{\left[\frac{n-1}{2}\right]}F_{l}(u^{2}+u^{2m})^{l}\right).$$
Suppose that $\Sigma_{1}=(0,+\infty)\backslash \{u\in(0,+\infty)|M_{0}(u)=0\}$. Let
$$M_{1}(u)=:\frac{M(u)}{M_{0}(u)},~~M_{1}^{'}(u):=\frac{M_{2}(u)}{M_{0}^{2}(u)},~u\in\Sigma_{1},
\eqno(4.9)$$
For $n<m$, by Lemma 4.1 we can obtain,
$$M_{2}(u)=\sum_{l=0}^{\left[\frac{n}{2}\right]+\left[\frac{n-1}{2}\right]}\sum_{j=0}^{l+1}\widetilde{C}_{l,j}u^{2l+2j(m-1)}+\sum_{l=0}^{2\left[\frac{n-1}{2}\right]}
\sum_{j=0}^{l+1}\widetilde{D}_{l,j}u^{2l+2j(m-1)+m},\eqno(4.10)$$
$\widetilde{C}_{l,j}$ and $\widetilde{D}_{l,j}$ are constants. It is obvious (4.10) hold when $2\left[\frac{n-1}{2}\right]<k-1$.

When $2\left[\frac{n-1}{2}\right]=k-1$, we have
$$\begin{aligned}
M_{2}(u)&=\sum_{l=0}^{\left[\frac{n}{2}\right]+\left[\frac{n-1}{2}\right]}\sum_{j=0}^{l+1}\hat{C}_{l,j}u^{2l+2j(m-1)}+\sum_{l=-\delta_{n}}^{\frac{m}{2}-2}
\sum_{j=0}^{l+1}\hat{D}_{l,j}u^{2l+2j(m-1)+m}\\
&+\sum_{j=1}^{\frac{m}{2}}\hat{D}_{\frac{m}{2}-1,j}u^{2(j+1)(m-1)}.
\end{aligned}\eqno(4.11)$$
When $2\left[\frac{n-1}{2}\right]>k-1$, we have
$$\begin{aligned}
M_{2}(u)&=\sum_{l=0}^{\left[\frac{n}{2}\right]+\left[\frac{n-1}{2}\right]}\sum_{j=0}^{l+1}\hat{C}_{l,j}u^{2l+2j(m-1)}+\sum_{l=\left[\frac{n}{2}\right]+\left[\frac{n-1}{2}\right]
-\frac{m}{2}+1}^{\frac{m}{2}-2}\sum_{j=0}^{l+1}\hat{D}_{l,j}u^{2l+2j(m-1)+m}\\
&+\sum_{l=\frac{m}{2}-1}^{2\left[\frac{n-1}{2}\right]}\sum_{j=l-\frac{m}{2}+2}^{l+1}\hat{D}_{l,j}
u^{2l+2j(m-1)+m},
\end{aligned}\eqno(4.12)$$
$\hat{C}_{l,j}$ and $\hat{D}_{l,j}$ are constants. Notice the order set $(u^{n_{1}},u^{n_{2}},u^{n_{3}},...,u^{n_{k}-1},u^{n_{k}})$ is an ECT-system on $u\in(0,+\infty)$ for $n_{i}\in\mathbb{N}$, hence we have
$$Z(m,n)\leq\begin{cases}
 4\left[\frac{n}{2}\right]^{2}+\left(6\delta_{n}+11\right)\left[\frac{n}{2}\right]+\frac{1}{2}\delta_{n}\left(5\delta_{n}+17\right)+4,~if~2\left[\frac{n-1}{2}\right]<\frac{m}{2}-1,\\\\
k^{2}+\frac{1}{2}\left((7-2\delta_{n})k-3-4\delta_{n}\right),~if ~2\left[\frac{n-1}{2}\right]=\frac{m}{2}-1,\\\\
      (4k+1)\left[\frac{n}{2}\right]+(3k+1)\delta_{n}-k(k-5)-1,~if ~2\left[\frac{n-1}{2}\right]>\frac{m}{2}-1,\\

\end{cases}$$
by using Rolle's Theorem. This ends the proof. $\diamondsuit$

\vskip 0.2 true cm
\noindent{\bf Proof of the case of} ${\bf (m,n)\in D_{5}\cup D_{6}}.$
\vskip 0.2 true cm

Similarly, the coefficients of $\bar{\Phi}_{k}$(reps. $\tilde{\Phi}_{k}$), $\bar{\Psi}_{k}$(reps. $\tilde{\Psi}_{k}$), $\bar{\Upsilon}_{k}$(reps. $\tilde{\Upsilon}_{k}$) and \begin{footnotesize}$${\bf \overline{G}_{2}}=:\frac{\partial{(\bar{\xi}_{0}^{(1)},\bar{\xi}_{1}^{(1)},...,\bar{\tau}_{0}^{(1)},\bar{\tau}_{1}^{(1)},...,\bar{\xi}_{0}^{(\frac{n-m}{2}+3)},\bar{\xi}_{1}^{(\frac{n-m}{2}+3)},...,\bar{\chi}_{1},\bar{\tau}_{0}^{(\frac{n-m}{2}+3)},...,
\bar{\eta}_{1},\bar{\xi}_{0}^{(\frac{m}{2}+2)},...,\bar{\eta}_{\frac{n-m}{2}+1})}}
{\partial{(\rho_{0,0},\zeta_{0,1},...,\gamma_{2t_{2}+1,0},\gamma_{2t_{2}+3,0},...,\gamma_{1,n-m+3},\gamma_{3,n-m+3},...,\rho_{n-m+2,1},\gamma_{1,n-m+4},...,\rho_{0,1},\gamma_{1,m+1},...,\rho_{n-m,1})}}
$$\end{footnotesize}

$$
  \addtocounter{MaxMatrixCols}{10}
    =\begin{pmatrix}
    \begin{smallmatrix}
          &2 &0 &\dots &0 &0 &\dots &0 &0 &\dots&0 &0 &\dots &0 &0 &\dots &0\\
          &0 &\pi\nu_{0,0}^{0} &\dots &0&0 &\dots &0 &0 &\dots &\frac{1}{2}\pi\nu_{0,0}^{0} &0 &\dots &0 &0 &\dots &0\\
          &\dots &\dots &\dots &\dots &\dots &\dots &\dots &\dots &\dots &\dots &\dots &\dots &\dots &\dots &\dots &\dots\\
          &0 &0 &\dots &-2 &0 &\dots &0 &0 &\dots &0 &0 &\dots &0 &0 &\dots &0\\
          &0 &0 &\dots &0 &-2 &\dots &0 &0 &\dots &0 &0 &\dots &0 &0 &\dots &0\\
          &\dots &\dots &\dots &\dots &\dots &\dots &\dots &\dots &\dots &\dots &\dots &\dots &\dots &\dots &\dots &\dots\\
          &0 &0 &\dots &0 &0 &\dots &-\frac{2}{n-m+4} &0 &\dots &0 &0 &\dots &0 &0 &\dots &0\\
          &0 &0 &\dots &0 &0 &\dots &0 &-\frac{2}{n-m+4} &\dots &0 &0 &\dots &0 &0 &\dots &0\\
          &\dots &\dots &\dots &\dots &\dots &\dots &\dots &\dots &\dots &\dots &\dots &\dots &\dots &\dots &\dots &\dots\\
          &0 &0 &\dots &0 &0 &\dots &0 &0 &\dots &2\nu_{\frac{n-m}{2}+1,0}^{0} &0 &\dots &0 &0 &\dots &0\\
          &0 &0 &\dots &0 &0 &\dots &0 &0 &\dots &0 &-\frac{2}{n-m+5} &\dots &0 &0 &\dots &0\\
          &\dots &\dots &\dots &\dots &\dots &\dots &\dots &\dots &\dots &\dots &\dots &\dots &\dots &\dots &\dots &\dots\\
          &0 &0 &\dots &0 &0 &\dots &0 &0 &\dots &0 &0 &\dots &2\nu_{0,0}^{0} &0 &\dots &0\\
          &0 &0 &\dots &0 &0 &\dots &0 &0 &\dots &0 &0 &\dots &0 &-\frac{2}{m+2} &\dots &0\\
          &\dots &\dots &\dots &\dots &\dots &\dots &\dots &\dots &\dots &\dots &\dots &\dots &\dots &\dots &\dots &\dots\\
          &0 &0 &\dots &0 &0 &\dots &0 &0 &\dots &0 &0 &\dots &0 &0 &\dots &2\nu_{\frac{n-m}{2},0}^{0}\\
       \end{smallmatrix}
       \end{pmatrix}
  $$

hence $\rm{det}{\bf \overline{G}_{2}}\neq 0$ since $\nu_{l,0}^{0}\neq 0$, which implies the independence of the coefficients.

For $n=2m-2,~2m=1$ and $2m$,
\begin{tiny}$${\bf \overline{G}_{3}}=:\frac{\partial{(\tilde{\xi}_{0}^{(1)},\tilde{\xi}_{1}^{(1)},...,\tilde{\tau}_{0}^{(1)},\tilde{\tau}_{1}^{(1)},...,\tilde{\eta}_{1},\tilde{\xi}_{0}^{(\frac{n-m}{2}+3)},...,\tilde{\chi}_{1},\tilde{\tau}_{0}^{(\frac{n-m}{2}+3)},...,
\tilde{\eta}_{\frac{m}{2}-1},\tilde{\xi}_{0}^{(\frac{n}{2}+1)},...,\tilde{\eta}_{\frac{m}{2}})}}
{\partial{(\rho_{0,0},\zeta_{0,1},...,\gamma_{2t_{2}+1,0},\gamma_{2t_{2}+3,0},...,\rho_{n-2m+2,1},\gamma_{1,n-m+3},...,\rho_{n-m+2,1},\gamma_{1,n-m+4},...,\rho_{n-m-2,1},\gamma_{1,n-1},...,\rho_{n-m,1})}}
$$\end{tiny}

$$
  \addtocounter{MaxMatrixCols}{10}
    =\begin{pmatrix}
    \begin{smallmatrix}
          &2 &0 &\dots &0 &0 &\dots &0 &0 &\dots&0 &0 &\dots &0 &0 &\dots &0\\
          &0 &\pi\nu_{0,0}^{0} &\dots &0&0 &\dots &0 &0 &\dots &0 &0 &\dots &0 &0 &\dots &0\\
          &\dots &\dots &\dots &\dots &\dots &\dots &\dots &\dots &\dots &\dots &\dots &\dots &\dots &\dots &\dots &\dots\\
          &0 &0 &\dots &-2 &0 &\dots &0 &0 &\dots &0 &0 &\dots &0 &0 &\dots &0\\
          &0 &0 &\dots &0 &-2 &\dots &0 &0 &\dots &0 &0 &\dots &0 &0 &\dots &0\\
          &\dots &\dots &\dots &\dots &\dots &\dots &\dots &\dots &\dots &\dots &\dots &\dots &\dots &\dots &\dots &\dots\\
          &0 &0 &\dots &0 &0 &\dots &2\nu_{\frac{n}{2}-m+1,0}^{0} &0 &\dots &0 &0 &\dots &0 &0 &\dots &0\\
          &0 &0 &\dots &0 &0 &\dots &0 &-\frac{2}{n-m+4} &\dots &0 &0 &\dots &0 &0 &\dots &0\\
          &\dots &\dots &\dots &\dots &\dots &\dots &\dots &\dots &\dots &\dots &\dots &\dots &\dots &\dots &\dots &\dots\\
          &0 &0 &\dots &0 &0 &\dots &0 &0 &\dots &2\nu_{\frac{n-m}{2}+1,0}^{0} &0 &\dots &0 &0 &\dots &0\\
          &0 &0 &\dots &0 &0 &\dots &0 &0 &\dots &0 &-\frac{2}{n-m+5} &\dots &0 &0 &\dots &0\\
          &\dots &\dots &\dots &\dots &\dots &\dots &\dots &\dots &\dots &\dots &\dots &\dots &\dots &\dots &\dots &\dots\\
          &0 &0 &\dots &0 &0 &\dots &0 &0 &\dots &0 &0 &\dots &2\nu_{\frac{n-m}{2}-1,0}^{0} &0 &\dots &0\\
          &0 &0 &\dots &0 &0 &\dots &0 &0 &\dots &0 &0 &\dots &0 &-\frac{2}{n} &\dots &0\\
          &\dots &\dots &\dots &\dots &\dots &\dots &\dots &\dots &\dots &\dots &\dots &\dots &\dots &\dots &\dots &\dots\\
          &0 &0 &\dots &0 &0 &\dots &0 &0 &\dots &0 &0 &\dots &0 &0 &\dots &2\nu_{\frac{n-m}{2},0}^{0}\\
       \end{smallmatrix}
       \end{pmatrix}
  $$

 hence $\rm{det}{\bf \overline{G}_{3}}\neq 0$ since $\nu_{l,d}^{0}\neq 0$, which implies the independence of the coefficients. Therefore, there exist $\rho_{i,j}$, $\gamma_{i,j}$ and $\zeta_{i,j}$ such that the lower bound of $Z(m,n)$ in Theorem 1.2 (ii) can be reached.

Now, we consider the upper bound for $(m,n)\in D_{5}\cup D_{6}$. For $m\leq n<2m$, we can obtain
$$\begin{aligned}
M(u)&=\sum_{l=0}^{\left[\frac{n}{2}\right]}\sum_{j=0}^{l}C_{l,j}u^{2l+2j(m-1)+1}+\sum_{l=\left[\frac{n}{2}\right]-\frac{m}{2}+1}^{\frac{m}{2}-1}\sum_{j=0}^{l}D_{l,j}u^{2l+2j(m-1)+m+1}\\
&+\sum_{l=\frac{m}{2}}^{\left[\frac{n-1}{2}\right]}\sum_{j=l-\frac{m}{2}+1}^{l}D_{l,j}u^{2l+2j(m-1)+m+1}+(u^{2}+u^{2m})\left(\sum_{l=0}^{\left[\frac{n-1}{2}\right]}E_{l}(u^{2}+u^{2m})^{l}\right)\\
&+(u^{2}+u^{2m})\left(\sum_{l=0}^{\left[\frac{n-1}{2}\right]}F_{l}(u^{2}+u^{2m})^{l}\right)\int_{0}^{1/\sqrt{1+u^{2m-2}}}\sqrt{1-t^{2}}dt.
\end{aligned}\eqno(4.13)$$
In this case,
$$\begin{aligned}
M_{2}(u)&=\sum_{l=0}^{\left[\frac{n}{2}\right]+\left[\frac{n-1}{2}\right]}\sum_{j=0}^{l+1}\widetilde{C}_{l,j}u^{2l+2j(m-1)}+\sum_{l=\frac{m}{2}}^{\left[\frac{n-1}{2}\right]+\frac{m}{2}}
\sum_{j=1}^{l+1}\widetilde{D}_{l,j}^{(1)}u^{2l+2j(m-1)+m}\\
&+\sum_{l=\left[\frac{n-1}{2}\right]+\frac{m}{2}+1}^{2\left[\frac{n-1}{2}\right]}\sum_{j=l-\left[\frac{n-1}{2}\right]-\frac{m}{2}+1}^{l+1}\widetilde{D}_{l,j}^{(1)}u^{2l+2j(m-1)+m}
+\sum_{l=\left[\frac{n}{2}\right]+\frac{m}{2}-1}^{\left[\frac{n-1}{2}\right]+\frac{m}{2}-1}\sum_{j=0}^{l+1}\widetilde{D}_{l,j}^{(2)}u^{2l+2j(m-1)+m}\\
&+\sum_{l=0}^{2\left[\frac{n-1}{2}\right]}\sum_{j=0}^{l}\widetilde{F}_{l,j}u^{2l+2j(m-1)+3m-2}\\
&=\sum_{l=0}^{m-2}\sum_{j=0}^{l+1}\hat{C}_{l,j}u^{2l+2j(m-1)}+\sum_{l=m-1}^{\left[\frac{n}{2}\right]+\left[\frac{n-1}{2}\right]}\sum_{j=l-m+2}^{l+1}\hat{C}_{l,j}u^{2l+2j(m-1)}\\
\end{aligned}$$
$$\begin{aligned}
&+\sum_{l=\left[\frac{n}{2}\right]+\left[\frac{n-1}{2}\right]-\frac{m}{2}+1}^{2\left[\frac{n-1}{2}\right]}\sum_{j=l-\frac{m}{2}+2}^{l+1}\hat{D}_{l,j}^{(1)}u^{2l+2j(m-1)+m},\quad\quad
\quad\quad\quad\quad\quad\quad\quad\quad\quad
\end{aligned}\eqno(4.14)$$
$\widetilde{C}_{l,j}$, $\widetilde{D}_{l,j}^{(i)}~(i=1,2)$, $\widetilde{F}_{l,j}$, $\hat{C}_{l,j}$ and $\hat{D}_{l,j}^{(1)}$ are constants.

For $n>2m$, we can obtain
$$\begin{aligned}
M(u)&=\sum_{l=0}^{m-1}\sum_{j=0}^{l}C_{l,j}u^{2l+2j(m-1)+1}+\sum_{l=m}^{\left[\frac{n}{2}\right]}\sum_{j=l-m+1}^{l}C_{l,j}u^{2l+2j(m-1)+1}\\
&+\sum_{l=\left[\frac{n}{2}\right]-\frac{m}{2}+1}^{\left[\frac{n-1}{2}\right]}\sum_{j=l-\frac{m}{2}+1}^{l}D_{l,j}u^{2l+2j(m-1)+m+1}+(u^{2}+u^{2m})\left(\sum_{l=0}^{\left[\frac{n-1}{2}\right]}
E_{l}(u^{2}+u^{2m})^{l}\right)\\
&+(u^{2}+u^{2m})\left(\sum_{l=0}^{\left[\frac{n-1}{2}\right]}F_{l}(u^{2}+u^{2m})^{l}\right)\int_{0}^{1/\sqrt{1+u^{2m-2}}}\sqrt{1-t^{2}}dt.
\end{aligned}\eqno(4.15)$$
In this case,
$$\begin{aligned}
M_{2}(u)&=\sum_{l=0}^{\left[\frac{n-1}{2}\right]+m-1}\sum_{j=0}^{l+1}\widetilde{C}_{l,j}^{(1)}u^{2l+2j(m-1)}+\sum_{l=m}^{\left[\frac{n-1}{2}\right]+m}\sum_{j=1}^{l+1}\widetilde{   C}_{l,j}^{(2)}u^{2l+2j(m-1)}\\
&+\sum_{l=\left[\frac{n-1}{2}\right]+m+1}^{\left[\frac{n}{2}\right]+\left[\frac{n-1}{2}\right]}\sum_{j=l-\left[\frac{n-1}{2}\right]-m+1}^{l+1}\widetilde{C}_{l,j}^{(2)}u^{2l+2j(m-1)}\\
&+\sum_{l=\left[\frac{n}{2}\right]+\left[\frac{n-1}{2}\right]-\frac{m}{2}+2}^{2\left[\frac{n-1}{2}\right]}\sum_{j=l-\left[\frac{n-1}{2}\right]-\frac{m}{2}+1}^{l+1}\widetilde{D}_{l,j}
u^{2l+2j(m-1)+m}+\sum_{l=0}^{2\left[\frac{n-1}{2}\right]}\sum_{j=0}^{l}\widetilde{F}_{l,j}u^{2l+2j(m-1)+3m-2}\\
&+\sum_{l=\left[\frac{n}{2}\right]-\frac{m}{2}+1}^{\left[\frac{n}{2}\right]+\left[\frac{n-1}{2}\right]-\frac{m}{2}+1}\sum_{j=\left[\frac{n}{2}\right]-m+2}^{l+1}\widetilde{D}_{l,j}
u^{2l+2j(m-1)+m}\\
&=\sum_{l=0}^{m-2}\sum_{j=0}^{l+1}\hat{C}_{l,j}^{(1)}u^{2l+2j(m-1)}+\sum_{l=m-1}^{\left[\frac{n-1}{2}\right]+m-1}\sum_{j=l-m+2}^{l+1}\hat{C}_{l,j}^{(1)}u^{2l+2j(m-1)}\\
&+\sum_{l=\left[\frac{n-1}{2}\right]+m}^{\left[\frac{n}{2}\right]+\left[\frac{n-1}{2}\right]}\sum_{j=l-m+2}^{l+1}\hat{C}_{l,j}^{(2)}u^{2l+2j(m-1)}+\sum_{l=\left[\frac{n}{2}\right]
+\left[\frac{n-1}{2}\right]+1}^{2\left[\frac{n-1}{2}\right]+\frac{m}{2}}\sum_{j=l-m+2}^{l-\frac{m}{2}+1}\hat{D}_{l,j}u^{2l+2j(m-1)},
\end{aligned}\eqno(4.16)$$
$\widetilde{C}_{l,j}^{(i)}~(i=1,2)$, $\widetilde{D}_{l,j}$, $\widetilde{F}_{l,j}$, $\hat{C}_{l,j}^{(i)}~(i=1,2)$ and $\hat{D}_{l,j}$ are constants. Thus we can get the results of Theorem 1.2 (ii) similar to the prove for case of $(m,n)\in D_{4}$. $\diamondsuit$

\end{document}